\documentstyle[amstex,12pt,epsf]{article}
\author{
D. V. Chistyakov \thanks
{ Kazan State University. E-mail: Dmitry.Chistyakov\symbol{64}ksu.ru. 
The work is partly supported by grants INTAS 00-334 and RFBR-01-02-17682-а.}
}

\title {
Fractal Measures,  $p$-Adic Numbers And
Continues Transition Between  Dimensions.
}

\renewcommand{\appendixname}{Appendix}
\date{}

\begin{document}
\hyphenation{trans-for-mat-ion}


\newcommand{\cA }{{\cal A}}
\newcommand{\cB }{{\cal B}}
\newcommand{\cC }{{\cal C}}
\newcommand{\cF }{{\cal F}}
\newcommand{\cE }{{\cal E}}
\newcommand{\cG }{{\cal G}}
\newcommand{\cD }{{\cal D}}
\newcommand{\cH }{{\cal H}}
\newcommand{\cK }{{\cal K}}
\newcommand{\cL }{{\cal L}}
\newcommand{\cU }{{\cal U}}
\newcommand{\cP }{{\cal P}}
\newcommand{\cO }{{\cal O}}
\newcommand{\cR }{{\cal R}}
\newcommand{\cS }{{\cal S}}
\newcommand{\cQ }{{\cal Q}}
\newcommand{\cJ }{{\cal J}}
\newcommand{\cI }{{\cal I}}
\newcommand{\cT }{{\cal T}}

\newcommand{\cV }{{\cal V}}
\newcommand{\cW }{{\cal W}}
\newcommand{\cY }{{\cal Y}}
\newcommand{\cX }{{\cal X}}
\newcommand{\cN }{{\cal N}}
\newcommand{\cZ }{{\cal Z}}
\newcommand{\cM }{{\cal M}}

\newcommand{\WL }{{\cA}}

\newcommand{\WI }{\WL (I)}

\newcommand{\fF}{f_{\cal F}}
\newcommand{\gG}{g_{\cal G}}


\def\figurename {{\rm \footnotesize Figure }}

\newtheorem{propdef}{Определение-предложение}
\newtheorem{mydef}{Определение}
\newtheorem{myprop}{Proposition}
\newtheorem{statement}{Proposition}
\newtheorem{mylemma}{Lemma}
\newtheorem{myso}{Corollary}
\newtheorem{mytheorem}{Теорема}
\newtheorem{mycond}{Соглашение}

\newcommand{\proof}{{\bf Proof~}}
\newcommand{\prooflem}{{\bf Proof \bf of Lemma.~}}
\newcommand{\cond}{{\bf Соглашение~}}
\newcommand{\conseq}{\\ \\{\bf Corollary~}}
\newcommand{\th }{{\rm th}}

\newcommand{\proofend}{$\Box$}
\newcommand{\Endproof}{\proofend\\}

\newcommand{\remark}{Замечание: }

\newcommand {\picr}[1]{ { \boxed {\rm #1}}}

\newcommand{\npage}{} 

\newcommand{\psfile}[2]{\centerline{\epsfxsize=#1 truecm\epsfbox{#2}}}

\newcommand{\bm}[1]{\{ #1 \} }

\newcommand{\Cond}[3]{\left \{ #1 \right \} ? \left [ #2 \right ] : \left [ #3  \right ] }

\newcommand{\rref}[1]{{\rm \ref{#1}}}
\newcommand{\bref}[1]{{\rm (}\rref{#1}{\rm)}}
\newcommand{\PS }[1]{{2^{#1}}}
\newcommand{\PX }{ \PS{X}}
\newcommand{\dW}{\ae}
\newcommand{\lW}{\ell}
\newcommand{\lWo}{\widetilde{\lW}}

\newcommand{\dist}[1]{ \nu[#1]}
\newcommand{\scr}[1]{ \sigma_{#1}}

\newcommand{\Diam}[1]{{\hat{#1}}}
\newcommand{\Drho}{\Diam{\rho}}

\newcommand{\ImU }{ \Upsilon}
\newcommand{\Ims }{ \ImU_s}

\newcommand{\Imm }{ \Upsilon_{s}^{(m)}}

\newcommand{\Imms }{ \widetilde{\Upsilon}_{s}^{(m)} }

\newcommand{\Immi }{ \Upsilon_{s}^{(\infty)}}
\newcommand{\Immn }[1]{ \Upsilon_{s}^{(#1)}}
\newcommand{\Immnn }[2]{ \Upsilon_{#2}^{(#1)}}

\newcommand{\wsol}{ \omega_{s}^{(m)}}
\newcommand{\wsoli}{ \omega_{s}^{(\infty)}}
\newcommand{\wsolip}{ \omega_{s/p}^{(\infty)}}
\newcommand{\Wsol}{ \Omega_{s,a}^{(m)}}
\newcommand{\Wsoli}{ \Omega_{s,a}^{(\infty)}}
\newcommand{\Wsolip}{ \Omega_{s/p,+0a}^{(\infty)}}
\newcommand{\Wsolipa}{ \Omega_{s/p,\ep   a}^{(\infty)}}
\newcommand{\idwsol}{ _{\otimes}\omega_{s}^{(m)}}
\newcommand{\facwsol}{  \mho_{s}^{(m)  }}
\newcommand{\wsoln }[1]{ \omega_{s}^{(#1)}}
\newcommand{\fp }[1]{\left \{#1\right \}_p}
\newcommand{\ip }[1]{\left [ #1\right ]_p}

\newcommand{\mK }[1]{k_{#1}}
\newcommand{\mKz }{ \mK{\zeta} }

\renewcommand{\cK }{\cC}

\newcommand{\hin }[1]{ \chi_{n}^{(m)}(#1)}
\newcommand{\hi }[1]{ \chi_{#1}^{(m)}}
\newcommand{\hins }{ \tilde \chi_{n}^{(m)} }
\newcommand{\hinn}[2]{\chi_{n}^{(#1)}(#2)}
\newcommand{\hinsn}[1]{\tilde \chi_{n}^{(#1)} }
\newcommand {\norm}[2]{\mid \mid #1 \mid\mid_{#2}}
\newcommand {\np}[1]{|#1|_{p}}
\newcommand {\npa}[1]{\mid #1 \mid_{p}^\alpha}
\newcommand {\abs}[1]{ |#1|}
\newcommand {\half}{ \frac{1}{2} }
\newcommand{\EL}{\stackrel{L}\simeq}
\newcommand{\Ll}{\stackrel{L}\preceq}
\newcommand{\LL}{\stackrel{L}\preceq}
\newcommand{\Lg}{\stackrel{L}\succeq}

\newcommand{\El}{\stackrel{\ell}\simeq}
\newcommand{\lL}{\stackrel{\ell}\preceq}

\newcommand{\lG}{\stackrel{\ell}\succeq}

\newcommand{\Eo}{\stackrel{\circ}\simeq}
\newcommand{\ol}{\stackrel{\circ}\preceq}
\newcommand{\oo}{\prec\prec}

\newcommand{\lel}{\lL}
\newcommand{\eql}{\El}
\newcommand{\leL}{\Ll}
\newcommand{\eqL}{\EL}
\newcommand{\leo}{\ol}
\newcommand{\eqo}{\Eo}
\newcommand{\leoo}{\oo}

\newcommand{\og}{\stackrel{\circ}\succeq}

\newcommand{\id }{{\rm id}}
\newcommand{\Ran }{{\rm Ran}}
\newcommand{\RE }[1]{{\rm Re}\left ( #1 \right )}
\newcommand{\IM }[1]{{\rm Im}\left ( #1 \right )}

\newcommand{\Rs}{{\bf R}}
\newcommand{\eRs}{\bar {\Rs}}
\newcommand{\Su}{{\bf S_1}}
\newcommand{\sC}{{\bf C}}
\newcommand{\Zs}{{\bf Z}}
\newcommand{\Qs}{{\bf Q}}

\newcommand{\Qqp}{Q^{(p)}}
\newcommand{\Ip}{ {\rm Z}(p^\infty)}

\newcommand{\eQs}{\bar {\Qs}}

\newcommand{\Cs}{ \sC}
\newcommand{\Ns}{{\bf N}}
\newcommand{\eNs}{\bar {\bf N}}

\newcommand{\Tn}[1]{{\bf T}_{#1}}
\newcommand{\Qn}[1]{{\bf Q}_{#1}}
\newcommand{\Zn}[1]{{\bf  Z}_{#1} }

\newcommand{\as}{{\bf a}}
\renewcommand{\as}{\{ a\}}

\newcommand{\Qa}{\Qn{\as}}
\newcommand{\Za}{\Zn{\as}}
\newcommand{\Ta}{\Tn{\as}}

\newcommand{\Qp}{\Qn{p}}
\newcommand{\Zp}{\Zn{p}}
\newcommand{\Tp}{\Tn{p}}

\newcommand{\Qpd}[1]{ \Qp ^{ ( #1) } }
\newcommand{\Qpdd}{\Qpd{D}}

\newcommand{\FRp}[1]{ \Rs_{+}^{#1} }
\newcommand{\pic}[2]{ fig. \ref{#1}.#2 }
\newcommand{\pf}[1]{ fig. \ref{fpi}.#1 }
\newcommand{\ps}[1]{ fig. \ref{spi}.#1 }

\newcommand{\UC}[1]{\rm U_{#1}}

\newcommand{\UCn}[1]{\UC{#1}^{\circ}}
\newcommand{\Ue}{\rm \UC{1}}
\newcommand{\Uen}{\rm \UCn{1}}

\newcommand{\diam}{\rm diam}
\newcommand{\la}{\lambda}
\newcommand{\ep}{\varepsilon}
\newcommand{\tilF}{\widetilde{F}}
\newcommand{\tilP}{\widetilde{P}}
\newcommand{\tilG}{\tilde G}
\newcommand{\tilC}{\widetilde{C}}
\newcommand{\tilT}{\widetilde{T}}
\newcommand{\tila}{{\tilde a}}
\newcommand{\tilf}{{\tilde f}}
\newcommand{\tilY}{\widetilde{Y}}
\newcommand{\tilZ}{\widetilde{Z}}
\newcommand{\tilO}{\widetilde{O}}
\newcommand{\tilPhi}{\widetilde{\Phi}}
\newcommand{\myind}{\mbox{\rm Ind\,}}
\newcommand{\FI}{$\,\mbox{Ф.И.}\!$}
\newcommand{\FIp}{\mbox{Ф.И.$p$-А}}
\newcommand{\tilR}{\widetilde{R}}
\newcommand{\calR}{{\cal R}}

\renewcommand{\Rs}{{\mathbb{R}}}
\renewcommand{\sC}{{\mathbb{C}}}
\renewcommand{\Zs}{{\mathbb{Z}}}
\renewcommand{\Qs}{{\mathbb{Q}}}


\maketitle
\begin{abstract}
\noindent
Fractal measures of images of continuous maps from the set of $p $-adic numbers $\Qp$
into complex plane $\sC$ are analyzed. Examples of "anomalous" fractals, i.e. the sets
where the $D$-dimensional Hausdorff measures (HM) are trivial, i.e. either zero, or
$\sigma$-infinite ($D $ is the Hausdorff dimension (HD) of this set) are presented.
Using the Caratheodory construction, the generalized scale-covariant HM
(GHM) being non-trivial on such fractals are constructed.
In particular, we present an example of $0$-fractal, the continuum
with HD$ =0 $ and nontrivial GHM invariant w.r.t. the group
of all diffeomorphisms $\sC$.
For conformal transformations of domains in $ \Rs^n $,
the formula for the change of variables for GHM  is obtained.
The family of continuous maps $ \Qp $ in $ \sC $
continuously dependent on "complex dimension"  $d \in \sC $ is obtained.
This family is such that:
1) if $d=2 (1) $, then the image of $\Qp$ is $\sC $ (real axis in $\sC $);
2) the fractal measures coincide with the images
of the Haar measure in $ \Qp $, and at $d=2 (1) $ they also
coincide with the flat (linear) Lebesgue measure;
3) integrals of entire functions over the fractal measures of images
for any compact set in $ \Qp $
are holomorphic in $d $, similarly to the dimensional regularization method
in QFT.
\end{abstract}

It is well-known that the Hausdorff measures (HM) are natural integral geometry characteristics
for a wide class of sets in $ \Rs^d $\cite {bill, FED, Feder}.
Therefore, contraction of a $D$-dimensional HM $h ^ {D} $ on  $D$-dimensional rectifiable submanifolds
is a measure of their areas \cite {FED} and, besides, there  exist fractal subsets
such that the HM contraction onto them for non-integer $D $ is also nontrivial,
i.e. is non-zero and ($ \sigma $-) finite, and determines their $D $-dimensional
fractal measures.
For each $F \subseteq \Rs^n $, there is a unique number $D _ {h} (F) $ called
Hausdorff dimension (HD) such that $h ^{D} (F) = \infty $ at $D < D _ {h} (F) $ and
$h ^ {D} (F) =0 $ at $D > D _ {h} (F) $.
However, HM contraction onto $F $ at $D=D _ {h} (F) $ can be trivial \cite {bill}.
Then, one may naturally ask
whether it is still possible in this anomalous case
to construct such a nontrivial measure $ \mu $ for $F $
which would possess the basic properties specific to the HM.
These are the following properties:
1) The measure must be obtained basically in the same way as
HM was, but with a generically wider class of test functions, i.e.
must be obtained by the  Caratheodory construction  \cite {FED};
2) The measure must not depend on any scale parameter but on the metric.
Thus, it follows from the dimensional analysis that the measure should be {\it scale covariant}:
$ \forall A \subset \Rs^n $, $ \forall \lambda > 0 $ and for some $D > 0 $
$ \mu (\lambda A) = \lambda^D \mu (A)$.
In this paper we construct such measures. They  have a unique jumping point
$D_F $ on each $F$ similarly to the HM, and always $D_F=D _ {h} (F)$.
If $D=0 $, then $h^D $ is a countable measure
that is nontrivial only on final and countable sets.
However, we show below that there exists a class of scale invariant measures ($D=0 $) such that
these measures are nontrivial on continuums and, moreover,
they are invariant with respect to any diffeomorphisms $ \Rs^n $.
It turns out that it is rather simple to gain and study examples of such sets by
considering them as the images of continuous embedding of  $p$-adic numbers field
$ \Qp $ in $ \Rs^n $ \cite {Zelenov, VVZ, Pit, Chis}.
This is basically due to ultrametricity of $ \Qp $ \cite {Kob, GG, VVZ} that
the $p $-adic counterparts of the real situation become much
simpler and sometimes even more correctly defined \cite {VVZ, Pit, MIS}.
Such spaces enjoy one important property:
any monotone increasing function of ultrametrics
is again an ultrametrics with equivalent uniform structure
(and, therefore, with equivalent topology) \cite {Kob}.
Thus, varying the metric, it is possible to assign to the same subset $F \subseteq \Qp $
the originally given Hausdorff dimension (HD) such that the fractal measure of $F$ coincides
with the contraction of the Haar measures on $\Qp$ \cite {Chis} onto $F$.
Therefore, if for some uniform embedding $ \ImU:\Qp \mapsto \Rs^n $
one selects a metric in $ \Qp $ such that it is in a sense close to the metrics
induced from $ \Rs^n $
(see section \ref {FPADIC}),
the examination of fractal properties of $ \ImU (F) $ can be reduced to examination of the properties of
$F \subseteq \Qp $.
In this paper, the uniform continuous maps $ \ImU_s:\Qp \mapsto \sC $, possessing the scaling property
$ \ImU_s (px) =s \ImU_s (x) $ with $s \in \sC $ are studied. Under some additional assumptions,  (analyticity in $s $, uniformity etc.)
they enjoy a series of remarkable properties.
It appears that any such $s$-parametric set of $\ImU_s $ is uniquely determined by the function
$\phi: [0,1] \mapsto \sC$ such that its continuity is sufficient for the fractal measure of the embedding image $ \ImU_s $ to coincide with the image of the Haar measure in $ \Qp $.
We construct examples of anomalous fractals which are images
of  multiple derivatives and integrals of $ \ImU_s $ over the parameter $s$.
The class of  $ 0 $-fractals being the images of embeddings of the corresponding continuum subsets
of $ \Qp $ with nontrivial 0-dimensional fractal measures is presented.
Some examples of the sets in $ \sC $ of zero  Lebesgue measure  having the HD$ =2 $ is also constructed.
A family of continuous maps from $ \Qp $ to $ \sC $ which
continuously depend on the parameter $d \in \sC $ is obtained.
This family is such that:
1) if $d=2 (1) $, then the image of $\Qp$ is $\sC $ (the real axis in $\sC $);
2) the fractal measures coincide with the images
of the Haar measure in $ \Qp $ and, at $d=2 (1) $, they also
coincide with the flat (linear) Lebesgue measures;
3) integrals of entire functions
over the fractal measures of  images of any compact set in $ \Qp $
are holomorphic in $d $.
Thus, the 1- and 2- dimensional integrals of holomorphic functions can be interpreted
as values of a function holomorphic in $d $, much similarly to the method of
dimensional regularization of  Feynman integrals \cite {Speer, Hooft, Wilson} \footnotemark.
In addition, note that the values of these functions at noninteger  $d $ are interpreted as integrals
over the corresponding fractal measures.
\footnotetext {
It is necessary to point out that the similarity between these concepts is rather formal.
}
\npage

\section{Pseudometric space $\WI$ }\label{FILTERS}

When examining properties of measures, metrics and test functions given on different spaces,
it is convenient to define two functions $f$ and $g$ to be $\varepsilon$-close if
$e^{-\varepsilon} f(x) < g(x) < e^{+\varepsilon}  f(x)$ for $x$ from a set
such that
for metrics it is "a set of infinitely close points" and for test
functions it is "an infinitely small set" only.
The construction considered below allows one to describe within a uniform scheme
the spaces of measures, test functions, metrics etc.
as quasiorder pseudometric spaces. It is important, however,
that the  Caratheodory construction is a functor, i.e. is a quasiorder-preserving contracting map
from the space of test functions into the space of exterior measures.
~

Let $I$ be any set, then, denote by $\WI$ the set of all pairs $\fF \equiv (f,\cF)$,
where  $f: I\mapsto \eRs^{+} $
\footnote{
$\eRs^{+}\equiv \eRs \cap \{x: x \ge 0 \}$, where $ \eRs \equiv \Rs \cup  \{  - \infty,\infty  \} $
is the extended number axis. It would be more naturally to define functions as
$f: V \mapsto \eRs^{+} $, where $V \in \cF$,
however, we can assume that  $f=\infty$ outside $V$.
From now on, we assume that
$\inf \{ \emptyset \}=  \infty$  and  $\sup\{ \emptyset \}=0 $ on $\eRs^{+}$
and
$\inf \{ \emptyset \}= \infty$  and  $\sup\{ \emptyset \}=- \infty$ on $\eRs$.
}
and $\cF$ is a filter at $I$.
\footnote{ The family $\cF $  of subsets of the set  $I$ is called
filter at $I$, if  $I\in \cF,\emptyset \notin  \cF $;
if $ A \in \cF$  and $ B \supseteq A$, then $B \in \cF$; $\forall A,B \in \cF,  A\cap B \in \cF$
\cite{Bur1} .}
Let us define the relations of quasiorder and a pseudometric at $\WI$ :
~

{\it
Relations $ \lL,\ol,\oo$ at $\WI$ :
$\fF \lL \gG$, $\fF \ol \gG$, $\fF \oo \gG $ if
$\lW(\fF,\gG)<\infty $,
$\lW(\fF,\gG) \leq 1$ and
$\lW(\fF,\gG)=0$,
respectively}, where $\forall \fF,\gG \in \WI$
\begin{eqnarray}
\lW (\fF,\gG)=
\left \{
\begin{array}{ccc}
\lWo (\fF,\gG) & \equiv \sup_{ U \in  \cG  } \inf_{ x \in U} \left ( r \in \Rs^{+} : f(x) \le r g(x))    \right )
&  ~if~  \cF \subseteq \cG  \\
\infty & ~~~ &  ~if~ \cF \not \subseteq \cG
\end{array}
\right. .
\label{LERE}
\end{eqnarray}
{ \it
Pseudometric $ \dW: \WI \times  \WI \mapsto  \eRs^{+}$:
$\dW(\fF,\gG) =   \ln \left ( \max\left ( \lW(\fF,\gG),\lW(\gG,\fF)    \right )   \right ) $.
}
We write {\it $\fF \El (\Eo) \gG$, if
$\fF \lL (\ol)\gG $  and  $\gG\lL (\ol)\fF $ at the same time }(cf. with $O$-symbolics).
It is obvious that $\fF \El (\Eo) \gG $, if and only if $\ae(\fF,\gG)<\infty (= 0)$\footnotemark.
\footnotetext{
The partition of $\WI$ according to the equivalence relation $\El$ is exactly the
partition of topological space $\WI$ into connected components. Indeed,
$\{h_{\cH} \in \WI : \dW(h_{\cH}, f_{\cF})<\infty \}$ ~$\forall \fF \in \WI $ is an open and close
set in $\WI$ at the same time,  and for any $f_{\cF}' \El \fF$
the map $[0,1] \owns t \mapsto (tf+(1-t)f',\cF) \in \WI$ is a continuous path
from $f_{\cF}$ to $f_{\cF}'$.}
It is easy to show that
\begin{eqnarray}
\ln \lWo(\fF,\gG)=
\inf \limits _{ \cG  \owns  U }  \sup\limits _{U \owns x  }   \{ \ln f(x)-\ln g(x) : f(x)\neq g(x) \},
\label{defLW}\\
\dW(f_\cG,\gG)=
\inf_{  \cG  \owns  U }  \sup \{ \abs{\ln f(x)-\ln g(x)} : x \in U\cap\{ f(x)\neq g(x)\}\},
\label{defDW}\\
\kappa(f_\cG,\gG)\equiv
\th \left ( \frac{\ae(f_\cG,\gG)}{2}  \right )=
{ \inf_{  \cG \owns  U }
\sup   \left \{   \left |  \frac{ f(x) - g(x) }{f(x) + g(x) } \right |: x \in U\cap\{ f(x)\neq g(x)\} \right \} }.
\label{defKap}
\end{eqnarray}
We shall consider only  {\it self-consistent subspaces in $\WI$},
i.e. such $\cS(I) \subseteq \WI$ that $\forall \fF,\gG \in \cS(I)$
$\lW(\fF,\gG)=\lWo(\fF,\gG)$ (i.e. if $ \lWo(\fF,\gG) <\infty $ then $ \cF \subseteq \cG $
and  $\fF \lL \gG$).
Let us consider two subsets in $\WI$. The first set is a subset of all pairs $(f,\cF)$ with
trivial $\cF$ : $\cF=\{I\}$. The second set $\cN(I)$ consists of $(f,\cF) \in \WL(I)$
such that $f^{-1}(\{0\}) \neq \emptyset$  and $\cF$ is the pre-image of the filter
of neighbourhoods of zero, which means that the set
$\{f^{-1}(U) : \forall~ open ~U :~ U\cap 0 \neq\emptyset  \}$
is the base of this filter
\footnote{The family of sets $\cB $ is called filter base if
$\emptyset \notin  \cB $;  $\forall A,B \in \cB,  \exists C \in \cB : C \subseteq A\cap
B$, and $\cF = \{B : \exists A \in \cB : A \subseteq  B    \}$ is called
filter generated by $\cB$.}.
It is easy to check that $\WL_0(I)$ and $\cN(I)$ are self-consistent.
Moreover, {\it $\WL_0(I)$ is a metric space}, the relations
$ \ol$ and $\Eo$  at $\WL_0(I)$ are the standard relations
$ \leq$ and  $=$, respectively, and if $f \oo g$, then either $f=0$, or $g=\infty$.
Let us introduce the following notations:
$\forall (f,\cF),(g,\cG) \in \WL(I)$  we write
$f \Ll (\EL) g$ if $(f,\{I\}) \lL (\El) (g,\{I\})$ or $\forall x_k \in I ~f(x_1,
\cdots x_k \cdots x_n) \leL(\eqL) g(x_1, \cdots x_k \cdots x_n)$, and
for any constant $c$ we write $c \eqL 1$
if $0<c<\infty$.

We call any map $\cC :\Sigma \mapsto  \WL(J')$ for  $\Sigma  \subseteq \WL(J)$
{\it  $\lL (\ol, \oo)$-isotonic}  map or {\it   $ \lL (\ol, \oo)$ -isotonia}, if
$\cC(\fF) \lL( \ol ,\oo) \cC(\gG)$  for all $\fF,\gG \in \Sigma$ such that
$\fF \lL( \ol ,\oo) \gG$, we also call just {\it isotonia} the $\lL,\ol,\oo,$-isotonic map.
For any set  $X$, denote by $\cD(X)$~$(\cD_0(X))$
the set of all $(\rho,\cF) \in  \cN(X \times X)~(\WL_0(X \times X))$
such that $\rho$ is pseudometric\footnotemark, and denote by
$\cM(X) \subseteq \WL_0(\PS{X})$ the set of all exterior
$\sigma$-semiadditive measures on $X$.
\footnotetext{
One can show that  $\cD_0(X)$
is a complete metric space.
}

Any map $\Phi : X \mapsto Y$ induces an isometry $\Phi^* : \cD(Y) \mapsto  \cD(X)$
such that
$ \forall \rho \in \cD(Y) ~\Phi^*(\rho)\equiv \rho_\Phi (\cdot,\cdot)= \rho(\Phi(\cdot),\Phi(\cdot))$.
Let us consider two maps
$\Phi_1 : X \mapsto (M_1,\rho_1) $ and  $\Phi_2 : X \mapsto (M_2,\rho_2) $,
where $M_i$  is metric spaces with metric $\rho_i$. Then, we define the
distance $\dW(\Phi_1,\Phi_2)$ between $\Phi_1 $ and $\Phi_2$ by the
formula
$$
\dW(\Phi_1,\Phi_2) \equiv \dW(\Phi_1^*(\rho_1),\Phi_2^*(\rho_2)),
$$
we also write $\Phi_1 \lL(\ol,\oo) \Phi_2$ as soon as $\Phi_1^*(\rho_1) \lL(\ol,\oo) \Phi_2^*(\rho_2)$.
We call a given map $\Phi : (M_1,\rho_1) \mapsto (M_2,\rho_2) $
{\it $ \ell$-contraction} if $\Phi \lL \id$, and {\it $\ell$-isometry}
if $\Phi \El \id$, where $\id$ is the identity map on $ M _ 1$.
It is clear that $\rho_2(\Phi(\cdot),\Phi(\cdot)) \lL (\eql) \rho_1(\cdot,\cdot)$,
iff $\Phi$ is  an $\ell$-contraction ($ \ell $ - isometry), in particular,
the $L$-contraction is the Lipschitzian map.
It is easy to see that the condition  $\Phi_1 \leo(\lel)  \Phi_2$ is equivalent
to the existence of the ($\ell$-) contraction $\Phi_{12} :  \Phi_2(X) \mapsto   \Phi_1(X)$,
such that $\Phi_{1}=\Phi_{12}\circ \Phi_{2}$
\footnote{It deserves noting that, if one considers a category $ {\cal K} (X)
$ such that all the objects in it are $ \Phi: X \mapsto (M, \rho) $ with
$ \Phi(X)= M $  and morphisms $ \lel $, the indicated correspondence is a functor from $ {\cal K} (X) $
into a category of metric spaces with uniform continuous maps between them.}.
Therefore, if $ \Phi _ 1 $ is injective, then $ \Phi _ 2 $ is also injective.
Besides, if there are given appropriate structures of either (uniform)
topological, or metric space,  on $ X $,
then, if $ \Phi _ 2 $ is a (uniform) continuous map
or $ \ell $ -contraction, the same is $ \Phi _ 1 $.

\section
{ The Caratheodory construction  }\label{KK}
~
Let  $\cJ $ be a family of subsets of any set  $X$.
For any pair $(\zeta ,\cF) \in  \WL(\cJ)$
we define an exterior $\sigma$-semiadditive measure $ \mKz$
by putting $ \forall A \subseteq X $,
\begin{eqnarray}
\mKz(A)=\sup_{ \cU \in \cF  }  \inf
\{   \sum_{i \in J} \zeta(S_{i}) :  S_i  \in \cU
, \{S_{i}\}_{i \in J} - {\it countable~cover~of~} A ~
\}.
\label{mhe}
\end{eqnarray}
This correspondence determines the contracting isotonia
$\cK  : \WL(\cJ) \mapsto \cM(X)$,
which we  call {\it Caratheodory construction} (CC) by analogy with the standard
Caratheodory construction  \cite{FED} which is a particular case of CC.

Hereafter, we always suppose that $ \cJ = \PX $.
Let us define the space of test functions $ \cT(X) \subseteq  \WL (\PX)$
as a set
of pairs $ (\zeta, \cF) $ such that the following natural requirements hold
\begin{eqnarray}
\forall U \in \cF  ~ if~ A \subseteq B \in U ~ then ~A \in U \label{test1}
\footnotemark \\
\exists V \in \cF  ~such\ that ~
if ~ A \subseteq B \in V ~then ~\zeta(A) \leq \zeta(B)\label{test2}
\end{eqnarray}
\footnotetext{Figuratively speaking, any subset of a "U-small" set is also
"U-small".}
It is clear that $\cM(X)  \subseteq \cT(X) $ and  $\cK_{|\cM(X)}=\id$.

With each pseudometric $\rho \in \cD(X)$
we associate a pair $(\Drho ,\cF_\rho) \in \cN(\PX) \subset \WL(\PX)$ such that
$\forall S \subseteq X~ \Drho(S)\equiv\diam(S)=\sup_{x,y\in S}(\rho(x,y)) $
and $\cF_\rho$ is the filter generated by the base
$\{ S\in \cJ  : \Drho(S)\leq 1/n\}_{n=1,2...}$.
Property (\ref{test1},\ref{test2}) for $(\Drho ,\cF_\rho)$ obviously
holds and
$\lW (\Diam{\alpha},\Diam{\beta}) \leq \lW ( \alpha ,\beta)$
$\forall \alpha,\beta \in \cD(X)$.
Therefore, the map $\Diam{~}: \cD(X) \mapsto  \cT(X) $ is an isotonic contraction.
Moreover, for any monotone
non-decreasing function $\zeta : \eRs^+ \mapsto \eRs^+$, using $\Diam{~}$
one can associate a map
$ \Diam{\zeta}: \cD(X) \mapsto  \cT(X) $ such that
\begin{eqnarray}
\Diam{\zeta}_{\rho}(\cdot )=(\zeta(\Diam{\rho}(\cdot) ),\cF_\rho) \label{testr}
\end{eqnarray}
It can be shown that $\mK{ \Diam{\zeta}_{\rho}}$ is
regular in the sense of Borel \cite{FED}.
Let $\zeta(r)=r^d$, then for each
pseudometric space $(X, \rho )$
the composition $\cK \circ \Diam{\zeta} :  \cD(X) \mapsto  \cM(X) $ defines
a $d$-dimensional
Hausdorff measure $h^d$ (HM)
(up to a constant factor $2^{-d}\Gamma(\frac{1}{2})^d  / \Gamma(1+\frac{d}{2})$
-- "the area of $d$--dimensional  sphere";
however, as far as we are interested in normalized fractal measures, this factor is
inessential).
For $D_1>D_2$,
since $\cK$ is an isotonia and $\Diam{\rho}^{D_1} \oo \Diam{\rho}^{D_2}$, we have
$h^{D_1} \oo h^{D_2}$. One obtains from the last inequality
the well-known formula/definition for the Hausdorff dimension (HD) \cite{bill,Feder} :
$ \forall A \subset X$
\begin{eqnarray}
D_{h}(A)\equiv \inf \{ \delta :  h^{\delta}(A)= 0 \} =\sup \{ \delta :  h^{\delta}(A)=  \infty \}  \label{HD}
\end{eqnarray}

For any map $\Phi : X \mapsto Y$ consider now the map $\Phi^* : \cT(Y) \mapsto  \cT(X)$
such that $\forall (\eta,\cF) \in \cT(Y)$  $\Phi^*\eta(S) = \eta(\Phi(S)) $, and
$ \Phi^*\cF$ is the filter generated by the base
$\{S \subseteq  X : \exists O \in \cU : S  \subseteq  \Phi^{-1}(O)   \}_{\cU \in \cF}$.
Let $\Phi : \cM(X) \mapsto  \cM(Y)$ be the map such that
$\forall \mu  \in \cM(X) ~\Phi(\mu)(S)=\mu( \Phi^{-1}(S))$, $\forall B \subset Y$
and $ \mu(S)|_{B}\equiv \mu(S \cap B)$. Then, from (\ref{test1},\ref{test2})
it follows that
the diagram
\begin{eqnarray} \label{diagr}
\end{eqnarray}
\begin{picture}(120,35)(-70,-55)
\put(43,-0) {\makebox(0,0){ $\Diam{\zeta}$ }}
\put(5,-10) {\makebox(0,0){ $\cD(Y)$ }}
\put(90,-10) {\makebox(0,0){$\cT(Y)$ }}
\put(25,-10){\vector(1,0){40}}
\put(7,-20){\vector(0,-1){22}}
\put(82,-20){\vector(0,-1){22}}
\put(0,-30){\makebox(0,0) {$\Phi^*$}}
\put(90,-30){\makebox(0,0) {$\Phi^*$}}
\put(43,-40) {\makebox(0,0){ $ \Diam{\zeta}$   }}
\put(10,-50) {\makebox(0,0){$\cD(X)$ }}
\put(90,-50) {\makebox(0,0){$\cT(X)$ }}
\put(25,-50){\vector(1,0){40}}
\put(185,-0) {\makebox(0,0){ $\cK$ }}
\put(255,-0) {\makebox(0,0){ $|_{ \Phi(X)}$ }}

\put(155,-10) {\makebox(0,0){ $\cT(Y)$ }}
\put(220,-10) {\makebox(0,0){$\cM(Y)$ }}
\put(175,-10){\vector(1,0){20}}

\put(290,-10) {\makebox(0,0){$\cM(Y)$ }}
\put(240,-10){\vector(1,0){30}}

\put(157,-20){\vector(0,-1){22}}
\put(290,-40){\vector(0,1){22}}
\put(150,-30){\makebox(0,0) {$\Phi^*$}}
\put(297,-30){\makebox(0,0) {$\Phi$}}
\put(223,-40) {\makebox(0,0){ $ \cK $   }}
\put(160,-50) {\makebox(0,0){$\cT(X)$ }}
\put(290,-50) {\makebox(0,0){$\cM(X)$ }}
\put(175,-50){\vector(1,0){90}}
\end{picture}

is commutative.
In particular, from \bref{diagr} it follows that the contraction of the CC measure
onto any $X \subseteq Y$ coincides with the CC measure on $X$ independently of $Y$.

First of all, we are interested in {\it scale-covariant measures}. These are the
measures $\mK{}$ such
that under the scale transform $ \lambda : X \mapsto X$  such
that $\rho( \lambda(x),  \lambda(y))=\lambda \rho(x,y)$ with some $\lambda>0$:
\begin{eqnarray}
\mK{}( \lambda(\cdot))=\lambda^D  \mK{}(\cdot) \label{mscaling}
\end{eqnarray}
for some {\it scale dimension} $D>0$.
It is clear that any HM  is a scale-covariant measure,
however, this is not the only possible choice.
Indeed, since $\cK$ is a contraction, then $\eta_{\cF} \Eo \zeta_{\cG}$ implies
that   $ \mK{ \eta_{\cF}}=\mK{\zeta_{\cG}}$. Therefore,
to satisfy condition \bref{mscaling} for $\mK{\zeta_{\cG}}$ it suffices
\begin{eqnarray}
\zeta_\cF( \lambda(\cdot))   \Eo  \lambda^D  \zeta_\cF(\cdot) \label{zscaling}.
\end{eqnarray}
For any $D \geq 0$ let us denote by $\cS^D $
{\it the set $\forall (\eta,\cN) \in \WL(\eRs)$ such that $\cN $ is the filter
of neighbourhoods of zero and  $\eta \eqo \varsigma$, where
$\varsigma$ is a monotone nondecreasing function such that
$\varsigma(r)>0$ with $r>0$ and}
\begin{eqnarray}\label{scalim}
\lim_{r \mapsto 0} \frac{\varsigma(\lambda r)}{\varsigma(r)}=\lambda^{D}.
\end{eqnarray}
The following lemma is correct (see proof in the Appendix):
\begin{mylemma} \label{isoton}
If $\zeta \in \cS^D$, then the
map $\zeta : \cN(I) \owns \fF  \mapsto (\zeta \circ  f,\cF) \in  \WL(I)$  at $D>0$ ($D=0$)
is isotonia  ( $ \lL$ ,$\ol$-isotonia) and  $\forall \fF, \gG \in \cN(I)$
if  $\lW(\fF, \gG)<\infty$, then
$$\lW(\zeta(\fF), \zeta (\gG))  \leq  \lW(\fF, \gG)^D .$$
In addition, if $D'>D$ and $ \fF \lel \gG$, then $\forall \eta \in \cS^{D'}$
$ \eta(\fF) \oo \zeta(\gG).$
\end{mylemma}

It is easy to show that $\forall \eta \in \cS^D$ and
$\forall \rho \in \cD(X)$, the function $\Diam{\eta_\rho}$
satisfies \bref{zscaling}.
Moreover, since   $\cK$  is an isotonia in accordance with lemma
\ref{isoton}, then it follows from lemma 1
\begin{statement}\label{DUni}
For $D\geq 0$ and  $\forall \zeta \in \cS^D$ we have\\
{\rm 1)} The map $\Diam{\zeta}:\cD(X) \mapsto \cT(X)$ (and also
$\cK \circ \Diam{\zeta}:\cD(X) \mapsto \cM(X)$) is $ \lL,\ol$-isotonic
and, if additionally  $D>0$, it is  $ \oo$-isotonic.\\
{\rm 2)} $\forall \alpha,\beta \in \cD(X)$
\begin{eqnarray}
\dW (\mK{\Diam{\zeta_\alpha}},\mK{\Diam{\zeta_\beta}}) \leq
\dW (\Diam{\zeta_\alpha},\Diam{\zeta_\beta}) \leq D \cdot \dW ( \alpha ,\beta).
\label{Lipsh}
\end{eqnarray}
{\rm 3)} $\forall \rho \in \cD(X)$ , $D_1,D_2\geq 0$ and
$\forall \zeta^1 \in \cS^{D_1}$, $\forall \zeta^2 \in \cS^{D_2}$,
if   $D_1>D_2$, then
\begin{eqnarray}
\Diam{\zeta}^{1}_{\rho}\oo\Diam{\zeta}^{2}_{\rho},~
\mK{\Diam{\zeta}^{1}_{\rho}}\oo \mK{ \Diam{\zeta}^{2}_{\rho}}.
\label{KDim}
\end{eqnarray}
\end{statement}

\begin{myso}\label{conseq1}
{\it If for some set $F$ $\exists \eta,\xi \in \cS^D$
such that $\mK{\eta}=0$ and  $\mK{\xi}>0$, then $D=D_h(F)$. In particular,
if $\exists \zeta \in \cS^D$ such that  $\mK{\zeta}$ is nontrivial on $F$,
then  $D=D_h(F)$. }
\end{myso}
\begin{myso}\label{conseq2}
Let  $(X,\rho)$ and $(Y,d)$ be metric spaces and   $\Phi : X \mapsto Y$.
Then, if $\Phi$ is an $\ell$-contraction ($\ell$-isometry), for any
$\zeta \in \cS^D$
all the arrows in    \bref{diagr}  are isotonic and Lipschitzian maps.
Therefore,
$$
\mK{\zeta_{d}}|_{\Phi(X)} \leL(\eqL)  \Phi \mK{\zeta_{\rho}},
$$
in particular, $\forall A \subseteq X$ ~$D_h(\Phi(A)) \leq (=) D_h(A)$.
In addition, if $D=0$, then $ \mK{\zeta_{d}}|_{\Phi(X)} \leq (=)  \Phi \mK{\zeta_{\rho}} $.
Therefore, if $\zeta_{0} \in \cS^{0}$, then $\mK{\zeta_{0}} $ is invariant
with respect to
all $\ell$-isometries, in particular, $\mK{\zeta_{0}} $ on $\Rs^n$ is invariant
w.r.t.
all diffeomorphisms of $\Rs^n$ (since $\Rs^n$ is a countable union of compact sets ) .
\end{myso}
In the next section, we give an example of $ \mK {\zeta _ {0}} $
which is non-trivial on a continuum.
The following claim also turns out to be valid:
\begin{statement}\label{Jac}
Let  $O,O'$ be  open domains  in   $\Rs^n$, $\Phi$ is
a conformal  map from $O'$ onto $O$,
${\rm J}(\Phi)$ denotes the Jacobian for  $\Phi$
and $\zeta_{D} \in \cS^D$.
Then, for any $\mK{\zeta_{D}}|_{O}$-summable function $f$,
for a change of variables the following formula holds:
\footnotemark
\begin{eqnarray}
\int_{O} f(y)
\mK{\zeta_{D}}(dy)=
\int_{O'} f(\Phi(x)) \left ( \sqrt[n]{\abs{{\rm J}(\Phi)(x)}} \right ) ^{D} \mK{\zeta_{D}}(dx)
\label{jacmes}.
\end{eqnarray}
\end{statement}
\footnotetext{
Note that, for $\Rs^2 \simeq \sC$, formula \bref{jacmes}
is valid for all biholomorphic $\Phi$, while on $\Rs^1$ it is correct
for all diffeomorphisms of $\Rs$. Moreover, in both cases
$\left ( \sqrt[n]{\abs{{\rm Jac}(\Phi)(x)}} \right ) ^{D}=\abs{\Phi'(x)}^D $.
}
\proof
For linear $\Phi$, the proof immediately follows from \bref{mscaling}.
In the generic case, the proof follows from statement 7 and since
due to lemma 3.2.2 \cite{FED} $\forall \delta>0$ there
exists a countable covering by Borel sets $ E_i$ such that
$\forall c \in E_i~ \dW(\Phi'(c)|_{E_i},\Phi|_{E_i}) < \delta$.
\Endproof
\npage
\section
{Fractal Measures   in  $\Qp$ .  }\label{FPADIC}
~
Fractals having an hierarchical structure
(for example, the Cantor set, the Sierpinski triangle, the Koch  curve etc.
\cite {Mandel, Feder})
are convenient to consider as images of uniformly continuous maps from
ultrametric spaces to $\Rs^d $.
Such fractals can be constructed through the following procedure.
Let $\cC_n$ be a sequence of finite or countable  families of compact
sets in $\Rs^d$, called clusters of level $n$, and
$\bm{a} \equiv\{a_n \}_{n\in \Zs}$ be a sequence such that
$\forall C_n \in \cC_n  ~ C_n =\bigcup \limits_{x_n=0}^{a_n-1} C_{n+1}^{x_n}$,
where $C_{n+1}^{x_n}  \in \cC_{n+1}$.
Let us assume that
$g_n \equiv \sup \limits_{C_n \in \cC_n} {\diam} (C_n) \to 0$ as $ n \to \infty$.
Let us fix any set $C_0 \in \cC_{-1}$, then the fractal required is the set
$\cF \equiv \bigcap \limits_{n \in \Ns} \{ \bigcup  {C: C \in \cC_n,C
\subseteq C_0} \}$.
Each point $f \in \cF$ is the limit of some thread
$C_0 \supseteq C_{1}^{x_0}   \supseteq C_{2}^{x_1} \supseteq .... C_{n+1}^{x_n} \supseteq...$
(i.e. $\bigcap\limits_{n=0}^{\infty} C_{n+1}^{x_n} = \{f\}$)
and, vice versa, the limit of each such thread  is an element of $\cF$.

Let us identify the set of such threads with the set $\Za$ of all sequences
$\{x_n \} _ {n=0} ^ {\infty}$ such that $x_n\in \{0..., a_n-1 \}$,
and also enter the (ultra)metric
$\rho (x, y) =g (v (x, y)) $ on $\Za$, where $v (x, y) = \sup \{n: x_n=y_n \}$
and {\it $g $ is any monotone decreasing function}, and $g (\infty) =0 $.
If $\ImU :\Za \mapsto \Rs^d$ is a map such that
$\bigcap \limits_{n=1}^{\infty} C_{n+1}^{x_n}=\{ \ImU(\{x_n\}_{n=0}^{\infty})
\}$,
then $\ImU$ is a uniformly continuous map on $\cF$, and if also
$\forall n ~g_n \Ll g(n)$, then $\ImU$ is the Lipschitzian map.
Let $ \ImU$ be an injective map  (for instance, this is the case, if the clusters
do not intersect).
If  $g (\cdot) $ is such that $ \ImU ^ {-1}:\cF \mapsto \Za $ is also
(locally) the Lipschitzian map,
then the construction of the fractal measure $ \cF $ can be transferred from $ \Rs^d $
onto the more convenient ultrametric space
$\Za $ isomorphic to the ring of $ \bm {a} $-adic integers \cite {HR}
and, in the special case of $a_n=p $, to the ring of $p$-adic integers
$\Zp \subseteq \Qp $.

Each element $x$ of the $p$-adic number field $\Qp$ is uniquely representable
as a formal power series  \cite{Kob},
\begin{eqnarray}\label{power}
x=\sum_{n=v}^{\infty}a_{n} p^{n}
= \sum_{n=v}^{-1}a_{n} p^{n}+\sum_{n=0}^{\infty}a_{n} p^{n}
\end{eqnarray}
with coefficients $a_{n} \in \{0,1,...,p-1\}$, where $v<\infty$ and
$p$ is some fixed prime number.
\footnote
{
Actually, the role of $p$ can be equally well played by any positive
integer,
since we nowhere use the existence of inverse elements in the ring $\Qp$.
}
The number $v(x)=v$ is called logarithmic norm of $x$.
Any  strictly monotone decreasing function $g$ such that $g(\infty)=0$
defines an invariant metric on the additive group $\Qp$ by the formula
$ \forall x,y \in  \Qp~\rho(x,y)\equiv\rho(x-y)  = g( v(x-y) )$
\footnote{
Note that the uniform structure (and, therefore, topology) in
$ \Qp $ does not depend on the choice of $g$.
}.
This metric  have {\it the ultrametric property}:
\begin{eqnarray}\label{ultra}
\rho (x-y)   \leq  \max( \rho (x) , \rho (y) ).
\end{eqnarray}
The series \bref{power}  absolutely converges in $\rho$.
Any number $q \in \Qs$ can be uniquely expanded into series \bref{power}
and $\Qp$ is the completion of $\Qs$ \cite{Kob}.
The first sum at the right-hand side of \bref{power} is denoted as $\fp{ x }$
being the fractional part of $x$. The second sum is denoted as $\ip{x}$
being the integer part of $x$. In this case,
$\fp{x} \in \Qs \cap [0,1) $ ~ and $\ip{x} \in \Zp$ ~,
where ~$ \Zp=\{ x \in \Qp:  \Vert x \Vert \leq 1  \} $~ is the ring of
$p$-adic integers
and $\np{x} \equiv p^{-v(x)}$ is {\it the canonical norm}.
We assume further (unless otherwise stated)
that there is a canonical norm in $ \Qp $.
It can be shown \cite{Chis} that {\it the 1-dimensional Hausdorff measure
in $(\Qp,\np{\cdot})$
coincides with the standard Haar measure $\chi$ in  $\Qp$}  such that
\begin{eqnarray}\label{nhaar}
\chi  (\Zp) =  \int  \limits_{ \Zp} d \chi = 1.
\end{eqnarray}
Let us consider $(\Qp,\rho)$ with $\rho(x,y)=g(v(x-y))$
and $\zeta(r) \eqo (\eql) p^{-g^{-1}(r)}$. Then, it is obvious that
\begin{eqnarray}\label{Taut}
\mK{\zeta_\rho}(\cdot)=(\eqL) \chi(\cdot).
\end{eqnarray}
Thus, $ \forall \zeta \in \cS^{D}$
one can find the proper metric in $ \Qp $; for instance, if
\begin{eqnarray}\label{TestFunc}
\zeta( r )  \eqo
\prod_{k=0}^{N}  \frac{1}{ \left (  \log_p^{\{ k \}}(1/r ) \right  )^{D_k}}
= \prod_{k=m}^{N}  \frac{1}{ \left (  \log_p^{\{ k \}}(1/r ) \right  )^{D_k}},
\end{eqnarray}
where $(\exp_a(t)\equiv a^t)$,
\begin{eqnarray}
\log_a^{\{ k \}}(t)\equiv\exp_a^{\{-k \}}(t)\equiv\left \{
\begin{array}{ccc}
\underbrace{\log_a(\log_a \cdots \log_a(t))\cdots))}_{k ~ times } &:& k>0 \\
t &:& k=0 \\
\underbrace{\exp_a(\exp_a \cdots \exp_a(t)\cdots))}_{-k~  times} &:& k<0\\
\end{array}
\right.
\end{eqnarray}
and $D_0=D$,  $D_m >0$ then
one can easily show that proper metric can be chosen as follows
\begin{eqnarray}\label{disteq}
\rho_{\{ \bar D\}}(x) \eqo
\frac{1}
{ \exp_p^{\{m\}}
\left (
\np{x}^{-\frac{1}{D_m}}
\prod \limits_{k=1}^{N-m}   \left ( \log_p^{\{ k \}}   \np{x}^{-1} \right )^{-\frac{D_{m+k}}{D_m}}
\right ) }.
\end{eqnarray}
Let us consider the map $J^m : \Qp \mapsto \Qp$ such that
\begin{eqnarray}\label{NULLMAP}
J^m \left (  \sum_{n=-\infty}^{\infty}x_{n} p^{n} \right )
= \sum_{n=-\infty}^{\infty} x_{n} p^{  \exp_p^{\{m\} }(n) }
\end{eqnarray}
is an $\ell$-isometry from $\Qp$, with the metric
$
\rho(x) =   p^{ - \exp_p^{\{m\}}\left ( v(x) \right ) }
$
in $(\Qp,\np{\cdot})$. Thus, for $m>0$ one obtains
a continuum $J^m (\Qp) \subseteq \Qp$
with HD$=0$ and the nontrivial measure $\mKz$
for $\zeta(r) \eqo \left ( \log_p^{\{m\}}({1}/{r}) \right )^{-1}$.
Moreover, $\mKz$ is invariant w.r.t. any $\ell$-isometry on  $\sC$.

So far all these constructions in $ \Qp $ look like a tautology,
however, they become nontrivial if one finds out
an $\ell$-isometry from $( \Qp, \rho _ {\{D \} }) $ to $ \Rs^d $.
For $ \Lambda_N\equiv \{x\in \Qp: \np {x} \leq p^N \} ~ (\Lambda_\infty =\Qp) $,
we call a continuous map $ \ImU: \Lambda_N \mapsto \Rs^d $ {\it automodel}
if for some $m\geq 0$
\begin{eqnarray}\label{selfsim}
\ImU(B_l^{n}) = \bigcup_{\bar l =0}^{p^m-1}
E_{l,\bar l,n} \circ \ImU(B_{\bar l}^{nm})  ,
\end{eqnarray}
$\forall n \in \Zs, l \in \Qp $. Here
$B_l^n\equiv \{x \in \Lambda_N: \np{x-l}\leq p^{-n}\}$ are
closed (and simultaneously open) balls in $\Lambda_N$
and $E_{l,\bar l,n} : \Rs^d\mapsto \Rs^d$ are isometry maps.
We call  $ \ImU: \Lambda_N \mapsto \Rs^d $  {\it quasiautomodel}
if there exists a sequence of automodel maps $\ImU_k$ such that
$\lim_{k \to \infty}\dW(\ImU_k,\ImU)=0$.
The following statement is correct:
\begin{statement}\label{AutoModel}
Let $\ImU: \Lambda_N \mapsto  \Rs^d$ be quasiautomodel,
$\zeta \in \cS^{D}$ is such that $\mKz(\ImU(\Zp)) < \infty$
and $\mKz(\ImU(B_l^n) \cap \ImU(B_{l'}^{n'}))=0$, when
$B_l^n \cap B_{l'}^{n'} =\emptyset$.
Then $\forall A \subseteq \Rs^d$
\begin{eqnarray}\label{KKHaar}
\mKz(A\cap \ImU( \Lambda_n) )=
\mKz|_{\ImU( \Lambda_n)}(A)=\mKz(\ImU(\Zp)) \chi(\ImU^{-1} (A) ),
\end{eqnarray}
\end{statement}
\proof
From the proposition  \bref{DUni} it follows that it is sufficient
to consider the automodel map.
For the automodel map, formula \bref {KKHaar} is valid
$\forall B_{\bar l}^{n} \subset \Qp$.
Indeed, the clusters $ \ImU (B_l ^ {n}) $ in \bref {selfsim} do not
overlap $\mKz $-nearly everywhere and $\mKz $ is a translation invariant.
Since $ \mKz$ and $ \chi$ are regular in Borel sense,
$\ImU (\Lambda_n)=\bigcup_{k<n} \ImU (\Lambda_k)$ is a Borel set
($\Lambda_k$ is a compact set, and $\ImU$ is a continuous map)
and the semi-ring $\{ B_l^{n} \}_{n,l \in \Zs}$ generates
a Borel $\sigma$-algebra in $\Lambda_N$, then
formula \bref {KKHaar} is valid $\forall A \subseteq \Rs^d$  .
\Endproof
Thus,  {\it if $\ImU$   is  quasiautomodel and
$\exists \zeta \in \cS^{D}$ such that  $\ImU$ is an $\ell$ -isometry
from $(\Qp,\rho )$ into $\Rs^d$
for some  $\rho(x,y)=g(v(x-y))$ such that  $\zeta(r) \eqo p^{-g^{-1}(r)}$,
then the fractal measure of $\ImU(\Lambda_N)$ defined by
\begin{eqnarray}\label{FrMes}
\mu_{\ImU(\Lambda_N)}(\cdot) = \frac{1}{\mKz(\ImU(\Zp))} \mKz( \Lambda_N \cap \cdot )
\end{eqnarray}
is the image of the Haar measure  in  $\Qp$ i.e.
$ \mu_{\ImU(\Lambda_N)}(\cdot)=\chi(\ImU^{-1} (\cdot) )$ }.
\npage
\section
{Maps $: \Qp \to \sC$.  }\label{PADIC}
~
The map $ \Qp \owns x \mapsto p x \in \Qp $ generates a
natural group of scaling transformations  $ \{p^n \} _ {n \in \Zs}(\sim \Zs)$
in the ring $ \Qp $   $ (\np {p^nx} =p ^ {-n} \np {x}) $.
On the other hand, any scaling transform of $ \sC $ is of the form
$ \sC \owns z \mapsto s z+t \in \sC $, where $s, t \in \sC $.
We want to describe {\it scaling-covariant} uniformly continuous maps
$ \ImU: \Qp \mapsto \sC$, i.e. those consistent with the scaling transforms:
$ \ImU (px) = s \ImU (x) +t ~ \forall x \in \Qp $ for some $s, t \in \sC $.
For the map $ \Ims (\cdot) = \ImU (\cdot) -\ImU (0) $ the latter requirement
reduces to the following principal condition on $ \Ims $
\begin{eqnarray}\label{scalling}
\Ims(px) = s \Ims(x) =p^{-\frac{1}{D_s}} e^{i \arg(s)}  \Ims(x), 
\end{eqnarray}
where $D_s=-\ln(p)/\ln\abs{s} $. Let us  assume that  $ \Ims (\cdot) \not =
0 $. Then from  \bref {scalling}, the continuity of $ \Ims $ and from the
condition $ \Ims (0)
= 0 $ it follows that $s \in \Ue ~ (\UC {r} \equiv \{ z \in \sC: \abs {z} <
r \}) $ \footnote{ Note that if $\Ims: \Qp \mapsto \sC$ is injective, then $
\lim_{\np {x} \rightarrow \infty} \abs {\Ims (x)} = \infty $. Therefore, one
could consider the continuous maps $ \ImU: \bar \Qp \mapsto \bar \sC $ such
that $ \ImU (px) = L (\ImU (x)) $, where $ \bar \Qp \equiv \Qp \cap \{\infty
\} $ is the single-point compactification of $ \Qp $, $ \bar \sC $ is the
Riemannian sphere and $ L (z) = (a z +b) / (c + d z) $. Then, $ \ImU (0) $
and $ \ImU (\infty) $ are fixed points of $L $, however, if $ \ImU (0)
\not = \ImU (\infty) $, there exists a linear-fractional automorphism
$U:\bar \sC \mapsto \bar \sC $ such that, for $ \Ims = U \circ \ImU $,
\bref {scalling} is valid \cite{Shabat}. }.

Let us define  two  numbers
\begin{eqnarray} \label{da}
\Delta^\pm (\Ims ) = \inf  \{ \abs{  \Ims(x) - \Ims(y) }^{\pm 1} :\forall x,y \in \Qp : \np{x-y}=1\}.
\end{eqnarray}
From \bref{scalling} it follows that
{\it  $\Delta^{-}(\Ims)>0$ iff   $\Ims$  is uniformly continuous
and $\Delta^{+}(\Ims)>0$ iff there exists a uniformly continuous map
$\Ims^{-1}: \Ims(\Qp) \mapsto \Qp$ such that $\Ims^{-1}\circ \Ims = \id $. }
Let $\rho_D(x,y)=\sqrt[D]{\np{x-y}}$ and $d_{\Ims}(x,y)= \abs{  \Ims(x) - \Ims(y) }$,
then it is easy to show  that
$\lW(\rho_{D_s},d_{\Ims})=1/\Delta^{+}(\Ims)$ and
$\lW(d_{\Ims},\rho_{D_s})=1/\Delta^{-}(\Ims)$ and,
hence,
$\dW(d_{\Ims},\rho_{D_s})=-\ln( \min(\Delta^{+}(\Ims),\Delta^{-}(\Ims)) )$.
Therefore, {\it $\Ims : \Qp \mapsto \sC$  is uniformly continuous
iff $\Ims : (\Qp,\rho_{D_s}) \mapsto \sC$ is an $\ell$-isometry.}
From this one immediately obtains that
{\it for any scaling-covariant uniformly continuous embedding
$\Ims : \Qp \mapsto \sC$ the Hausdorff measure $h^{D_s}$
is nontrivial on $\Ims(\Qp)$
and $D_h(\Ims(B))=D_s$ for  any open set $B \subset \Qp$ }.
Besides, it is clear that if  $\Ims $  is only a uniformly continuous map,
then  $D_h(\Ims(\Qp)) \leq D_s$.
Let us assume furthermore that, for each fixed $x \in \Qp $ ,~ $\ImU _ {(\cdot)} (x) $
is a function holomorphic on $ \Uen ~ (\UCn {r} \equiv \UC {r} \backslash \{0 \})
$. Then the requirement $ \bref {scalling} $ and the residue at zero
uniquely determines the set of maps $ \{\Ims \} _ {s \in \Uen} $.
Indeed, let $ \phi (x) = \mathrel {\mathop {{\rm res}} _ {s=0}} (\Ims (x)) $, then
$\forall x \in \Qp$
\begin{eqnarray}\label{Loran}
\Ims(x) =\sum_{n=-\infty}^{\infty}  \phi  \left (\frac{x}{p^{n+1}}  \right   ) s^n.
\end{eqnarray}
Let us also assume that, for some $r > 0 $, there exists  $ \Delta_r > 0 $ such
that $ \Delta ^ {-} (\Ims) \geq \Delta_r ~ \forall s \in \UCn {r}$
\footnotemark.
Then, from the Cauchy inequalities for the Laurent series \bref{Loran}
it follows that $\phi(x+\Zp)=\phi(x)~ \forall x \in \Qp $,
or, which is the same, $\phi(x)=\phi(\fp{x})$.
\footnotetext{
One can easily show that the requirement
$ \Delta ^ {-} (\Ims) \geq \Delta_r > 0 $
is equivalent (up to the replacement $ \Ims \to s^M \Ims $)
to the equipotential uniform continuity of the set of functions $ \{\Ims \} _ {s
\in \UCn {r}} $.
}
Now the Laurent expansion of $ \Ims $ can be written as
\begin{eqnarray}\label{Loran2}
\Ims(x)=\Ims^{\phi}(x)=\sum_{n=v(x)}^{\infty}  \phi  \left ( \frac{x}{p^{n+1}} \right   ) s^n
=\sum_{n=-\infty}^{\infty}  \phi  \left ( \fp{\frac{x}{p^{n+1}}}  \right   ) s^n.
\end{eqnarray}
The function $ \phi $ is uniquely determined by its values on the set
$$\Ip\equiv \{ \fp{x}  : x \in \Qp \}=\{ l/p^k  : l,k \in \Ns, l<p^k \}.$$
Therefore, the function  $ \phi $  can be  considered as being
defined on $[0,1) \subset \Rs$ or
as a periodic function (with period $=1$ ) defined on $\Qs$, $\Qp$, or $\Rs$.
For the sake of simplicity, let us assume that
$ \phi: \Ip \mapsto \sC $ is an arbitrary bounded function
(without loss of generality, let $ \abs {\phi (x)} \leq 1 $).
Then, the series \bref{Loran2} defines the
map $\Ims^{\phi}$ such that \bref{scalling} holds and
$\Delta^{-}(\Ims^{\phi})\geq 1-\abs{s}$.
Using that $ \phi (x/p ^ {n+1}) $ depends only on $x_n, x _ {n-1}..., x _ {v (x)} $,
one easily gets the  following estimate for $\Delta^{+}(\Ims^{\phi})$:
\begin{eqnarray} \label{neqda}
\Delta^{+}(\Ims^{\phi}) \geq  \dist{\phi} -\frac{\abs{s}}{1-\abs{s}},
\end{eqnarray}
where
$\dist{\phi}=\inf_{a=1,p-1}
\left \{  \abs{\phi(p^{-1}a+ \tau)-\phi(\tau) }: \tau \in \Ip  \right  \}$
\footnote{
Considering $ \phi (t) $ as the periodic function $ \phi (t) = \phi (t+n) $,
it is possible to interpret $ \phi (t)$ as a trajectory of a particle in $ \sC $.
Then, $ \dist {\phi} $ is the minimal distance between $p $ particles
which started to move at sequent moments of "time"
$t=0, \frac {1} {p}..., \frac {p-1} {p} $.
}.
Thus, if $ \dist {\phi} > 0 $, then, for small enough $s $ such that
$ \abs {s} < \scr {\phi} \equiv {\dist {\phi}} / ({1 +\dist {\phi}}) $,
the map $ \Ims ^ {\phi} $ is an $\ell$-isometry.
For example, if $\phi(x)=(p-1)^{-1} x_{-1}$,
then  $\dist{\phi}=(p-1)^{-1}$ , $\scr {\phi}=p^{-1}$,
and if $\phi(x)=\exp(i 2 \pi \fp{x} )-1$,
then $\dist{\phi}=\sin(\pi/p)$ ,$\scr
{\phi}=(1+\sin(\pi/p)^{-1})^{-1} \geq p^{-1}$
(cf. with \cite{Chis}).

If $ \phi (x) = \phi (x _ {-1}) $, it is possible to obtain
a lower bound for  $D_h (\Ims ^ {\phi} (\Qp))$.
Indeed, it is easy to prove that $\ImU_{s^K}^{\phi}=\Ims^{\phi}\circ \Theta_K$,
where $\Theta_K:\Qp \mapsto \Qp$ such that
$\Theta_K\left (  \sum_{n=-\infty}^{\infty}x_{n} p^{n} \right )
= \sum_{n=-\infty}^{\infty} x_{n} p^{Kn}$.
Thus, for $K$ such that $\abs{s}^K<\scr{\phi}$, one gets
$$D_s/K=D_{s^K}=D_h(\Ims^{\phi}(\Theta(\Qp)))\leq D_h(\Ims^{\phi}(\Qp)).$$
For any integer $l $, let us introduce $ \forall x \in \Zp $ the map
$ \partial_s ^ {l} \Ims ^ {\phi}:\Zp \mapsto \sC $,
\begin{eqnarray}
\partial_s^{l}\Ims^{\phi}(x)=\sum_{n=l}^{\infty}  \phi  \left ( \fp{\frac{x}{p^{n+1}}} \right   )
\frac{n!}{(n-l)!} s^{n-l}.
\end{eqnarray}
Using \bref{app1},
it is easy to show that  $ d_{\partial_s^{l}\Ims^{\phi}}  \eql \rho_{ (D_s,l) }$
with   $\abs{s} < \scr{\phi}$,
where $\rho_{ (D_s,l)(x) }= v(x)^l p^{-v(x)/D_s }$ and
$d_{\partial_s^{l}\Ims}(x,y)= \abs{ \partial_s^{l} \Ims(x) -\partial_s^{l}\Ims(y)
}$.
Thus, {\it $\mKz$ with $\zeta(r)= (r\abs{\ln(r)}^{-l})^{D_s}$ ($\zeta \in \cS^{D_s}$)
is a nontrivial measure on $\partial_s^{l}\Ims^{\phi}(\Zp)$,
however, $h^{D_s}_{|\partial_s^{l}\Ims^{\phi}(\Zp)}=0$ at $l<0$
and is $\sigma$-infinite at $l>0$}.

Now let us construct {\it a set in $\sC$ of zero Lebesgue measure,
but with  HD$ =2 $}.
Put $\cF'_2=\partial_s^{-1}\Ims^{\phi} (\Zp)$  with $\phi(x)=x_{-1}$ and
$s=\frac{i}{\sqrt{p}}$.
Using \bref{app2}, one can prove that
$$\rho_{ (2,-1) } \lel  d_{\partial_s^{-1}\Ims^{\phi}}  \lel \rho_{(2,-2) }.$$
Hence, using corollary 1 of proposition \bref{DUni}, one obtains that
$D_h(\cF'_2)=2$ and $h^2(\cF'_2)=0$, since $\rho_{ (2,0) } \oo \rho_{ (2,-1) }$.

Since $\partial_s^{0}\Ims^{\phi}=\Ims^{\phi}|_{\Zp}$,
let us identify $ \partial_s ^ {0} \Ims ^ {\phi} $ with $ \Ims ^ {\phi}: \Qp \mapsto \sC$.
Let $[\phi]_m(x)\equiv \phi(p^{-m} \ip{p^m x})=\phi(\sum_{n=-m}^{0}x_{n} p^{n})$.
Then, similarly to \cite {Chis} it can be shown
that  $ \partial_s ^ {l} \Ims ^ {[\phi] _m} $ is automodel.
For any function  $\phi :[0,1) \mapsto \sC$, let
$\norm{\phi}{\infty}^{(p)} \equiv \sup_{} \{ \abs{\phi(q)}~ : q \in \Ip \} $.
A function  $\phi :[0,1) \mapsto \sC$ is called {\it $(p)$-continuous},
if $\lim_{m\to\infty}\norm{\phi-[\phi]_m}{\infty}^{(p)}=0$.
Obviously, $\forall M \in \Zs_{+}$
$\phi(x)=[\phi]_M(x)$  is a $(p)$-continuous function.
The following claim is correct.
\begin{statement}\label{SelfCont}
If $\partial_s^{l} \Ims^{\phi}$ is an $\ell$-isometry with
a $(p)$-continuous function $\phi$, then
$\partial_s^{l} \Ims^{\phi}$ is a quasi-automodel map
and $\exists M$ such that   $\partial_s^{l} \Ims^{[\phi]_m}$
is also an $\ell$-isometry at $m\geq M$.
\end{statement}
\proof
Since $ \dW (\rho _ {(D_s, l)}, d _ {\partial_s ^ {l} \Ims ^ {\phi}}) < \infty $
and for any $m \in \Zs _ {+} $ $\partial_s ^ {l} \Ims ^ {[\phi] _m}  $ is an
automodel map,
the proof follows from the following lemma (proved in the Appendix) applied to
$\phi_1=\phi$  and $\phi_2=[\phi]_m$.
\proofend
\begin{mylemma} \label{kappa}
For any bounded function $\phi_i: \Ip \mapsto \sC~(i=1,2)$
\begin{eqnarray}\label{divmap}
\kappa \left( \partial_s^{l} \Ims^{\phi_1},\partial_s^{l} \Ims^{\phi_2}\right )  \leq
c~  e^{\min(\dW_1,\dW_2)}  \norm{\phi_1-\phi_2}{\infty}^{(p)},
\end{eqnarray}
where
$\dW_i=\dW(\rho_{ (D_s,l)}, d_{\partial_s^{l}\Ims^{\phi_i}})$
and   $c={2}{(1-\abs{s})^{-1}}$.
\end{mylemma}
The requirement of $ (p)$-continuity is not too restrictive.
Indeed, the following assertion (proved in the Appendix) is valid.
\begin{statement}\label{pCont}
The function $ \phi: [0,1] \mapsto \sC $ is $ (p) $-continuous,
if it is right-continuous and is continuous everywhere except
for a finite number of points of the first kind discontinuity  $a_1..., a_N $
such that $a_i \in \Ip $, in particular, if $ \phi $ is just continuous
(continuity is understood in the sense of ordinary topology on $ \Rs $).
\end{statement}
From \ref{AutoModel} and \ref{SelfCont} it follows that
{\it $h^{D_s}(\Zp) \eqL 1$  at $ \Delta^{\phi}_s >0 $ and
the fractal measure $\Ims^{\phi}(\Qp)$ is the image of the Haar measure on $\Qp$.}
Similarly, the fractal measure $\partial_s^{l}\Ims^{\phi}(\Zp)$
(with $\zeta(r)= r^{D_s} \abs{\ln(r)}^{- D_s l} $ )
is the image of the Haar measure on $\Qp$ too.

Let now $ \Delta^{\phi}_s >0 $. It is easy to show that the map
$\cN \equiv \Ims^{\phi}  \circ J^m : \Qp \mapsto \sC $
is an $\ell$-isometry from $\Qp$ with metric
$ \rho_m (x) =   p^{ - \exp_p^{\{ m\}}\left ( \frac{v(x)}{D_s}\right ) }$
to $\sC$.
Let $\zeta(r) \eqo \left ( \log_p^{\{m\}}({1}/{r}) \right )^{-1} $, then $\zeta \in \cS^0$
and, from proposition \ref{DUni}, it follows that
$\mKz|_{\cN(\Qp) }(\cdot)=\mK{\zeta_{\rho_m} }|(\cN^{-1} (\cdot))=\chi(\cN^{-1} (\cdot))$
\footnote{
Note that  in this case
quasiautomodelity of $\cN$ is not used.
Moreover, $\mKz(\cN(\Zp))=\chi(\Zp)=1$
}.
The measure  $\mKz$ is invariant w.r.t.
any diffeomorphism $\sC \approx \Rs^2$,
and, at $m>1$, w.r.t. any homeomorphism with arbitrary Holder index.

Thus, in the cases considered above for any
$ \mu _ {\Ims ^ {\phi} (\Lambda)}$-summable function
$f: \sC \mapsto \sC $ one has
\begin{eqnarray}\label{FRACINT}
\int \limits_{\sC}   \mu_{\Ims^{\phi}(\Lambda)}(dz)~f(z)=
\int \limits_{\Lambda}  \chi(dx) f(\Ims^{\phi}(x))  .
\end{eqnarray}
Note that the integral in the right hand side of \bref {FRACINT} is correctly defined even
if the measure $ \mu _ {\Ims ^ {\phi} (\Lambda)} $ is not defined,
furthermore, the following proposition  holds (see proof in the Appendix).
\begin{statement}\label{DiffMes}
Let  $O $ be an open  set  in  $\UCn{r}$  for some $r<1$ ,
$\Lambda $  be a measurable subset in $ \Qp$,
$\ImU_{O}(\Lambda)\equiv \bigcup_{s \in O} \Ims(\Lambda) $  and
$f \in C^{m}(\ImU_{O}(\Lambda))$ and $l ,\bar l$ be such that
$l + \bar l \leq m$.
Then, if one of the following condition holds,
\\
{\rm 1)} the set $ \Lambda$ is bounded;\\
{\rm 2)}  $\delta_s^{\phi}\equiv \inf_{\np{x}=1}\abs{\Ims^{\phi}(x)}>0$
and $\exists \nu>0 $  such that with $ k\leq l, \bar  k \leq  \bar l $
the following inequality  holds
( $\partial_z \equiv \frac{\partial}{\partial x}+i \frac{\partial}{\partial y}$~
$\partial_{\bar z} \equiv \frac{\partial}{\partial x}-i \frac{\partial}{\partial y}$ )
$$
\forall z \in \ImU_{O}(\Lambda)~
\abs{\partial_z^k  \partial_{\bar z}^{\bar k} f(z)}
\leL \frac{1}{1+\abs{z}^{D_s+k + \bar k +\nu}};
$$
then the integral
\begin{eqnarray}\label{INT}
{\rm I}_{s}^{\Lambda}(f)\equiv \int \limits_{\Lambda} \chi(dx) f(\Ims^{\phi}(x)),
\end{eqnarray}
exists, ${\rm I}_{s}^{\Lambda}(f) \in C^{m}(O)$  and
the following equation is correct
\begin{eqnarray}\label{INTDIFF}
\partial_s^l  \partial_{\bar s}^{\bar l}{\rm I}_{s}^{\Lambda}(f)=
\int \limits_{\Lambda} \chi(dx)   \partial_s^l  \partial_{\bar s}^{\bar l} f(\Ims^{\phi}(x)).
\end{eqnarray}
In particular, if $f \in \cS(\Rs^2)$, where   $\cS(\Rs^2)$ is the Schwarz space
of rapidly decreasing  $C^{\infty}$-smooth functions,
then, ${\rm I}_{s}^{\Lambda}(f) \in C^{\infty}(O)$.
Furthermore, $\ImU_{O}(\Lambda)$ is an open set, since
from holomorphy of  $f$ in $\ImU_{O}(\Lambda)$
it follows that ${\rm I}_{s}^{\Lambda}(f)$ is
some holomorphic function in  $O$
$(\partial_{\bar z} f(z) =0 ~ \Rightarrow
\partial_{\bar s}{\rm I}_{s}^{\Lambda}(f)=0)$.
\end{statement}
The constructions considered in this section can be immediately generalized
onto maps from $ \Za $ to $ \sC $.
Let us define the map $\omega^{\phi}_\nu: \Rs \times \Za \mapsto \sC $
(assuming that  $\phi(t+1)=\phi(t)$) as follows
$\forall (\tau,x) \in \Rs \times \Za $
$$
\omega^{\phi}_\nu (\tau,x)=\sum_{n=0}^{\infty}
\left( \frac{1}{a^{(n)}} \right)^{\nu} \phi \left(  \frac{(x)_{a}^{n}+\tau}{a^{(n)}}            \right),
$$
where $a^{(n)}=\prod \limits_{k=0}^{n} a_k ~(a^{(-1)}\equiv 1)$
and  $  (x)_{a}^{n}=\sum \limits_{k=0}^{n} x_k a^{(k-1)}$.
In the special case of $a_k=p$, one has $\omega^{\phi}_\nu (0,x)=\Ims^{\phi}(x)$
$\forall x \in\Zp$ for $s=p^{-\nu}$.
Furthermore, it can be proved that, for example, if
$a_k=(k+2)!$, $\phi(t)=\exp(i2\pi t)$ and
$D^{-1}=\RE{\nu}>1$, then
$ \rho_{D,0}\lel d_{\tau} \lel \rho_{D,1}$,
where $\rho_{D,i}(x,y)=v(x,y)^i /\sqrt[D]{a^{(v(x,y))} }$
and $d_{\tau}(x,y)=\abs{\omega^{\phi}_\nu (\tau,x)-\omega^{\phi}_\nu (\tau,y)}$.
Therefore, $\forall \tau \in \Rs$ $D_h(\omega^{\phi}_\nu (\tau, \Za))=D$.
It is easy to show that
$\forall n\in\Zs ~ \omega^{\phi}_\nu (\tau,x)=\omega^{\phi}_\nu (\tau-n,x+n)$,
i.e. the function $ \omega ^ {\phi} _ \nu $ is constant on the cosets of the subgroup
$B\equiv\{(-n,n) \in  \Rs \times \Za : n\in\Zs \}$ and,
therefore, it can be represented as the map of
{\it $\bm{a}$-adic solenoid} \cite{HR}
$\Sigma_{\bm{a}}\equiv ( \Rs \times \Za)/B$ to $\sC$.
Similarly, the map $ \Omega ^ {\phi} _ {\nu, \alpha}: \Rs \times\Za \mapsto \Rs^3 $
such that $ \rho+i h =\omega ^ {\phi} _ \nu (\tau, x) + \alpha $ and
$ \phi =\tau $, where $h, \rho, \phi $ are cylindrical coordinates on $ \Rs^3
$, can be considered as the map from $ \Sigma _ {\bm {a}} $ into $ \Rs^3 $ \footnotemark.
\footnotetext{
Note that   $ \Sigma _ {\bm {a}} $ for $a_k = (k+2)! $ is isomorphic
to the character group of the additive group of the field $ \Qs $ \cite {HR}
or to the group $ {\bf A} / \Qs $, where $ {\bf A} $ is the ring of
adeles \cite {GG, Mam}.
}
If now $ \omega ^ {\phi} _ \nu $ is such that
$ \rho _ {D, 0} \lel d _ {\tau} \lel \rho _ {D, 1} $,
repeating the arguments from \cite {Chis} it can be shown
that, for $ \abs {\RE {\alpha}} $ large enough,
the map $ \Omega ^ {\phi} _ {\nu, \alpha} $ is a continuous embedding
and $D_h(\Omega^{\phi}_{\nu,\alpha}(\Sigma_{\bm{a}}))=D+1$.
Moreover, if  $\phi: [0,1] \mapsto \sC$ is a continuously differentiable
function, then $\forall (\tau,x) \in \Rs \times \Za $ the following equality holds
$$
\partial_\tau \omega^{\phi}_{\nu} (\tau,x)=  \omega^{\partial_\tau  \phi }_{\nu+1} (\tau,x)
$$
It follows that $\Omega^{\phi}_{\nu,\alpha}(\Sigma_{\bm{a}})$
is the invariant set of the autonomous ordinary differential equation in $ \Rs^3 $
with the Lipschitzian right hand side \cite {Chis}.
\npage
\section
{Continuous Transition between Dimensions}\label{TRANSIT}
Before considering a concrete set of maps realizing the transition between
integer dimensions
note  that the set $ \cF_s =\Ims (\Qp) $ is invariant w.r.t. the scale
transformation
$z \to s z $ and its fractal measure $ \mu _ {\cF_s} $ is transformed as follows
\begin{eqnarray}\label{Cdim}
\mu_{\cF_s}(s \cdot) =
{p}^{-1}\mu_{\cF_s}(\cdot) =
\abs{s^d}\mu_{\cF_s}(\cdot)
=e^{i\theta} s^d \mu_{\cF_s}(\cdot).
\end{eqnarray}
Here we introduce a "complex  dimension" $d \in \sC$ of the set
$\cF_{d}\equiv\cF_{s}$
and an arbitrary parameter $\theta \in \Rs$, with the following relations between
them being correct
\begin{eqnarray}
s=s(d) \equiv  \exp\left  (- \frac{\ln p+i \theta}{d} \right ) \label{Cdim01}.\\
D_s^{-1}=\RE{d^{-1}\left (1+i \frac{\theta}{\ln p} \right )} \geq 0  .
\end{eqnarray}

\begin{figure}
\psfile{14.0}{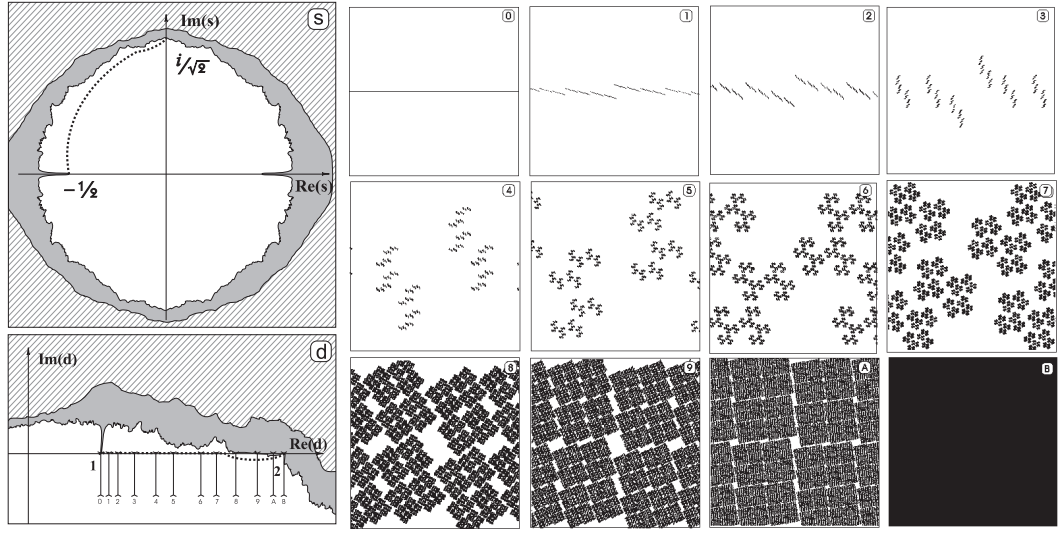}
\caption{
{\footnotesize
{ 
We depict by the dotted line the path from $d=1$ to $d=2$ in the $d$-plane
(fig. $ \picr {d} $)
and in the  $s$-plane (fig. $ \picr {s} $) for $p=2 $,   by white color
the  domain  of embeddings ($ \Delta ^ {-} (\Ims ^ {\phi}) > 2 ^ {-12} $)
and by grey color the domain  ($ \delta (\Ims ^ {\phi}) > 2 ^ {-10}) $.
In frames  {\footnotesize $ \picr {0} - \picr {B} $} the sets $ \cF_d =\ImU _ {s (d)} ^ {\phi} (\Qp) $
are shown for values $d $ lying on paths.
These values are marked  on a fig. $ \picr {d} $ by labels $ {\rm 0,..., B}$.
}
}}\label{splane}
\end{figure}


Let us consider now the map $ \Ims^{\phi}$ with $\phi$ such that
$\phi(x)=x_{-1}$ $\forall x=\sum \limits_{n=-\infty}^{\infty}x_{n} p^{n} \in \Qp$ .
Using the p-adic decomposition of the real numbers, it can be proved that
$\ImU_{\frac{1}{p}}^{\phi}(\Zp)=\{x+iy\in \sC :0 \leq x \leq \frac{p}{p-1}, y=0 \}$.
Now from \bref{scalling} it follows that
$\ImU_{\frac{1}{p}}^{\phi}(\Qp)=\{z \in \sC : \IM{z}=0 ,\RE{z}\geq 0 \}$.
It is easy to show that
$\ImU_{-\frac{1}{p}}^{\phi}(\Zp)-[\frac{p}{2}]
\sum \limits_{n=0}^\infty \left (\frac{-1}{p} \right)^n =
\ImU_{\frac{1}{p}}^{\phi}(\Zp)-[\frac{p}{2}] \sum_{n=0}^\infty
\left (\frac{1}{p} \right)^n$,  therefore,
$\ImU_{-\frac{1}{p}}^{\phi}(\Zp)=
\{x+iy\in \sC :0 \leq  x+\frac{2p^2 [p/2]}{p^2-1}  \leq \frac{p}{p-1}, y=0 \}$.
From this, using  \bref{scalling}, one gets
$\cF_d=\ImU_{-\frac{1}{p}}^{\phi}(\Qp)=\{z \in \sC : \IM{z}=0 \}$
with $d=1 \left( \frac{\ln p +i \theta}{\ln p +i {\pi}} \right) ~ (s(d)=- 1/p)$.
Let us define the bijection  $q: \Qp \times \Qp \mapsto \Qp$ such that
$$q(\sum_{n=-\infty}^{\infty} x_{n}  p^n,\sum_{n=-\infty}^{\infty} y_{n} p^n)
=\sum_{n=-\infty}^{\infty} x_{n} \left(p^{2}\right)^n
+p \sum_{n=-\infty}^{\infty} y_{n}\left(p^{2}\right)^n,
$$

One easily shows that $\ImU_{s}^{\phi}(q(x,y))=\ImU_{s^2}^{\phi}(x)+s\ImU_{s^2}^{\phi}(y)$
and, in particular,
$\ImU_{\pm \frac{i}{\sqrt{p}} }^{\phi}(q(x,y))=
\ImU_{-\frac{1}{p}}^{\phi}(x)  \pm i
\left( \frac{1}{\sqrt{p}} \ImU_{-\frac{1}{p}}^{\phi}(y) \right)$ .
Therefore, for
$d=2 \left( \frac{\ln p +i \theta}{\ln p \pm i {\pi}} \right)$
$ ( s(d)= \pm  \frac{i}{\sqrt{p}})$
one obtains that $\Ims^{\phi}(\Zp)$
is a closed rectangle in $\sC$ of the size
$ \frac{p}{p-1} \times  \frac{\sqrt{p}}{p-1}$
and $\Ims^{\phi}(\Qp)=\sC$.

Let us assume that $ \theta =\pi $, then one gets
that $\cF_1$ is the real axis  $ \sim \Rs^1$  and  $\cF_2 =\sC \sim \Rs^2$.
Besides, in spite of the fact that the map $ \ImU _ {s} ^ {\phi} $ is
not injective at $d=1,2 $, still $D_s=d $,
$h ^ {D_s} (\ImU _ {s (d)} ^ {\phi} (\Zp)) \eqL 1 $ and
the formula \bref {FRACINT} (with $ \zeta (\rho) = \rho ^ {D_s} $) is valid.
Thus, the fractal measures of  $\cF_1$ and $\cF_2$ are
(up to a constant factor) the 1-dimensional
$\mu_{\cF_1}(dz) \sim  \delta(z - \bar z) dz d\bar z$
and 2-dimensional $\mu_{\cF_1}(dz) \sim dz d\bar z$
Lebesgue measures, respectively.
It is easy to show that the points $ 1 $ and $ 2 $ in the $d$-plane can be
connected
by a continuous path such that $ \mu _ {\cF_d} (f) $ is defined on it
and $ \forall f\in\cS (\Rs^2)$,
due to proposition  \ref{DiffMes},
$\mu_{\cF_d}(f)$  $C^{\infty}$ is
a $C^{\infty}$ -smooth function and
if $f$ is a holomorphic function on $\UC{\sqrt{p}/(\sqrt{p}-1)}$,
then, $\mu_{\cF_d^0}(f)$ is a holomorphic function on this path,
$\cF_d^0=\ImU_{s(d)}^{\phi}(\Zp)$.
In figure \ref {splane} an example of such a path is shown in $d$- and $s$- planes
with $p=2$.
There are also drawn the (white) domain  of embeddings
$ \Delta ^ {-} (\Ims ^ {\phi}) > 2 ^ {-12} \approx 0 $,
where the fractal measures are defined and the formula \bref {FRACINT}  holds,
and the (grey) domain $ \delta (\Ims ^ {\phi}) > 2 ^ {-10} \approx 0 $,
where the requirements of assertion \ref {DiffMes} hold.

\npage
\section{\appendixname }\label{APPENDIX}
\newcommand{\Ts}{{\bf D}}
\proof {\bf of lemma \ref{isoton}}  ~
The following lemma holds
\begin{mylemma} \label{psiass}
If the function  $\phi : \Rs \mapsto \Rs$ such that
$\forall \delta \in \Rs$
$\lim \limits_{t \to \infty} \abs{\phi(t+\delta)-\phi(t)} = 0$
and $ \phi(t) = \mu (t) + \upsilon(t)$, where $\mu(t)$ and $\upsilon(t)$ are,
respectively, monotone and uniformly continuous functions
at some interval $(\tau_0, \infty )$,
then \\
\mbox{~~~}{\rm I)}~  $\forall \varepsilon>0 ~\exists \tau$ such that
$\abs{\phi(t)-\phi(t')}<\varepsilon +\varepsilon  \abs{t-t'}$ at
$t,t' \geq \tau$.\\
\mbox{~~}{\rm II)} ~$\forall r<\infty$
$
\lim \limits_{\tau \to \infty}\sup \{ \abs{\phi(t)-\phi(t')} : t,t'>\tau ,\abs{t-t'}< r\}  =0.
$\\
\mbox{~}{\rm III)}~  $\forall \varepsilon>0 ~\exists \tau$ such that
$\abs{\phi(t)}< \varepsilon t$ at $t \geq \tau$.
\end{mylemma}
\proof
$\forall \delta,\varepsilon>0 ~\exists \tau_{\delta,\varepsilon}>\tau_0$ such that
$\forall t' \geq \tau_{\delta,\varepsilon}$ $\abs{\phi(t'+\delta)-\phi(t')}< \varepsilon \delta/2$.
Let us choose $ \delta < 1 $ such that $\forall \delta' \leq \delta,\delta'>0$
~$\abs{\upsilon(t'+\delta')-\upsilon(t')} < \varepsilon/4$ .
Then, from monotonicity of $ \mu (t) $ it follows that
$\abs{\phi(t'+\delta')-\phi(t')} \leq \abs{\phi(t'+\delta)-\phi(t')} +\varepsilon/2 < \varepsilon$.
($
\abs{\phi(t'+\delta')-\phi(t')}-\varepsilon/4
\leq \abs{\mu(t'+\delta')-\mu(t')}
\leq \abs{\mu(t'+\delta)-\mu(t)}
\leq \abs{\phi(t'+\delta)-\phi(t')} +\varepsilon/4
$).
Let  $t\geq t'$ and $t_k=t +\delta k$, then $t-t_{N_{\delta}}=\delta'\leq \delta$,
where $N_{\delta}= \left[   {\delta^{-1}} ({t-t'})     \right]$
and the following chain of inequalities holds
\begin{eqnarray*}
\abs{\phi(t)-\phi(t')} \leq
\abs{\phi(t)-\phi(t_{N_{\delta}})}
+\abs{\phi(t_{N_{\delta}})-\phi(t')}\leq  \\
\leq    \abs{\phi(t_{N_{\delta}}+\delta')-\phi(t_{N_{\delta}})}
+\sum_{k=0}^{N_{\delta}-1} \abs{\phi(t_k)-\phi(t')}< 
\varepsilon +  \frac{ \varepsilon }{2} \delta N_{\delta}
\leq   \varepsilon + \frac{ \varepsilon }{2}  \abs{t-t'}.
\end{eqnarray*}
Now assertions I and II take place at  $\tau \geq \tau_{\delta,\varepsilon}$,
and assertion III can be obtained for
$\tau \geq   \abs{t'}+2+2\abs{\phi(t')}/ \varepsilon$.
\Endproof
Let us return now to the proof of lemma \ref {isoton}.
$\forall \fF, \gG \in \cN(I) ~ \ln \lW(\fF, \gG) = \ell_0 < \infty$,
Let $\alpha(s)=-\ln f(s)$ and $\beta(s)=-\ln g(s)$,
then  $\forall \epsilon>0~ \exists t_0$ such that $\forall s \in V_{t_0}=\{ s:  \beta (s) \geq t_0  \}$
$\delta(s)= \beta(s) - \alpha(s) < \ell=\ell_0+\epsilon$.
$\mu(t)=Dt+\phi(t)=-\ln\zeta(e^{-t})$
is a monotone increasing function, thus, $ \phi (t) $ obeys lemma  \ref{psiass}
and  $\exists t'_0>t_0$ such that
$\forall s \in V_{t'_0}~\delta\mu(s) = \mu(\beta(s)) -
\mu(\alpha(s))=D\delta(s)+\phi(\beta(s)) - \phi(\alpha(s))
\leq D\delta(s)+\epsilon \abs{\delta(s)} +\epsilon \leq
(D+\epsilon {\rm sign}{(\ell)})\ell +\epsilon $.
Therefore,  $\lW(\zeta(\fF), \zeta (\gG))  \leq  \lW(\fF, \gG)^D$ at $D>0$.
Since $\delta\mu(s)\leq 0$ at $\delta(s)\leq 0 $, then
$\delta\mu(s)\leq \max(\epsilon \delta(s) +\epsilon,0)$ at $D=0$.
From this one obtains from this
$\lW(\zeta(\fF), \zeta (\gG))\leq1$ at $\lW(\fF,\gG) <\infty$.
Let $D'>D$, then  $\forall \eta \in \cS^{D'}$
such that $\forall t>t_0$
$\ln \lW(\eta(\gG),\zeta (\gG)) \leq \ln(\eta(e^{-t})/\zeta(e^{-t})) =
-(D'-D)t +\phi(t)-\phi'(t)$.
Using assertion III of lemma \ref{psiass}, one gets that
$\ln \lW(\eta(\gG),\zeta (\gG))=-\infty$.
Thus,  $ \eta(\fF) \lel \eta(\gG) \oo \zeta(\gG)$.
\Endproof
Let $0\leq r < 1 $ and $p\in \Ns $,
then, using formula 5.2.3.1 \cite {PBM},
it can be proved  that the following relations hold at $v \to \infty $
\begin{eqnarray}
\sum_{n=v}^{\infty}\frac{n!}{(n-l)!} r^{n} \eqo
\frac{1}{1-r}
v^{-l} r^{v}
\label{app1}\\
\frac{ p^{-v}}{v^{2}}  \lel
\left(
\frac{ p^{-v}}{v}  -(p-1) \sum_{n=v+1}^{\infty} \frac{p^{-n} }{n}
\right)
\label{app2}.
\end{eqnarray}
\proof {\bf of lemma  \ref{kappa}}  ~
At $\dW_1=\dW_2=\infty$ the proof is trivial. Let now $\dW_1<\infty$.
Then $\cF_{\rho_{ (D_s,l)}}=\cF_{d_1}$, where
$d_i=d_{\partial_s^{l}\Ims^{\phi_i}}$.
$\forall \varepsilon>0~ \exists v_0$ such that
$  d_{i}(x,y) \geq e^{-(\dW_i+\varepsilon)} \rho_{ (D_s,l)}(x-y) $
at $ v(x,y)\geq v_0 $.
On the other hand,
$\Ims^{\phi_i}(x)-\Ims^{\phi_i}(y)=R_v^i(x)-R_v^i(y)$ at $ v(x,y)=v $,
where $R_v^i(x)=\sum_{n=v}^{\infty}  \phi_i  \left ( \fp{\frac{x}{p^{n+1}}} \right   )
\frac{n!}{(n-l)!} s^{n-l}$.
Therefore,
$\abs{d_1(x,y)-d_2(x,y)}\leq \abs{R_v^1(x)-R_v^2(x)}+\abs{R_v^1(y)-R_v^2(y)}$.
Using \bref{app1} and the equation $\rho_{ (D_s,l)}(x)=v(x)^l \abs{s}^{v(x)}$,
one can choose $v_0$ such that
$\abs{d_1(x,y)-d_2(x,y)} \leq e^{\varepsilon}
2/(1+\abs{s}) \norm{\phi_1-\phi_2}{\infty}^{(p)}  \rho_{ (D_s,l)}(x-y) $
at  $v>v_0$. Thus,
$\abs{d_1(x,y)-d_2(x,y)}/(d_1(x,y)+d_2(x,y))
\leq 2 (1-\abs(s))  e^{\min(\dW_1,\dW_2)+2\varepsilon}
\norm{\phi_1-\phi_2}{\infty}^{(p)}$.
Now the proof is obtained immediately from  \bref{defKap}.
\Endproof
\proof {\bf of proposition \ref{pCont}}  ~
If $ \phi $ is continuous on $[0,1]$,
then it is a uniformly continuous function and the proof is trivial.
Each $a_i $ has the form $a_i=l_i/p ^ {k (i)} $.
Let $m_0 >k(i)$, $a_0=0$, $a_{N+1}=1$
and $S_i=[a_{i},a_{i+1})$. It is easy to show that,
if $ q \in S_i \cap \Ip$, then also
$[q]^m\equiv p^{-m} \ip{p^m q}\in S_i$ at  $m>m_0$,
and, therefore, $\norm{\phi-[\phi]_m}{\infty}^{(p)} = \max_{i=0,N}
(\sup_{} \{ \abs{\phi(q)-[\phi]_m(q)}~ : q \in S_i \cap \Ip  \} ) $.
The function  $\phi$ is  uniformly continuous on  $S_i$,
since $\phi$ is continuous on $S_i$ and $\exists \lim_{x\to a_{i+1}-0} \phi(x)$.
Therefore, $\lim_{m\to\infty}\sup_{} \{ \abs{\phi(x)-[\phi]_m(x)}~ : x \in S_i \}=0$.
\Endproof
\proof {\bf of proposition \ref{DiffMes}}  ~
From theorem IV.115 \cite{Schwartz} it follows that, if there exists
a $\chi$-summable function $g$ such that $\forall x \in \Lambda $
$  \abs{\partial_s^l  \partial_{\bar s}^{\bar l} f(\Ims^{\phi}(x))} \leL
g(x)$, then there exists the integral ${\rm I}_{s}^{\Lambda}(f)$
and the equation  \bref{INTDIFF} holds.
One can easily show that, for
$\partial_s^l  \partial_{\bar s}^{\bar l} f(\Ims^{\phi}(x))$,
the following representation takes place
\begin{eqnarray} \label{LA1}
\partial_s^l  \partial_{\bar s}^{\bar l} f(\Ims^{\phi}(x))
=\sum \limits_{q,\bar q}^{l,\bar l} ( \partial_z^q  \partial_{\bar z}^{\bar q} f)(\Ims^{\phi}(x))
P_s^q(x) \bar P_s^{\bar q}(x),
\end{eqnarray}
where $P_s^q (x) $ are degree $q$ polynomials of
$ \partial_s ^ {k} \Ims ^ {\phi} (x) ~ (k=0..., q) $.
It can be proved that $\forall  q=0,1,2...$ the following inequalities hold
\begin{eqnarray} \label{LA2}
\sup_{\np{x}\leq 1}\abs{ {\partial_s^{q}\Ims^{\phi}} (x)} \leq \frac{ q!}{(1-\abs{s})^q}.
\end{eqnarray}
From this one obtains from these
$\sup_{x \in \Lambda_n}\abs{ \partial_s^{q}\Ims^{\phi} (x)} \leL  \abs{s}^{-n}$
$\forall n>0 $.
Therefore, if  $ \Lambda $ is bounded, then $ \ImU _ {O} (\Lambda) $ is bounded
too.
Now it follows from \bref {LA1} that
$\abs {\partial_s^l \partial _ {\bar s} ^ {\bar l} f (\Ims ^ {\phi} (x))} \leL 1 $
$ \forall x \in \Lambda $.
To prove assertion 2) note that
$
{\rm I}_{s}^{\Lambda}(f)={\rm I}_{s}^{\Lambda \cap \Zp}(f)
+\sum_{n=1}^{\infty}  {\rm I}_{s}^{\Lambda \cap {\rm S_n}}(f),
$
where ${\rm S_n} =\{x\in \Qp : \np{x}=p^n \}$.
Since $\abs{ \Ims^{\phi} (x)}\geq \delta_s^\phi \abs{s}^{-n}, ~  \forall x \in {\rm S_n}$,
one can prove using $\chi({\rm S_n})=\frac{(p-1)}{p}p^{n} $ that
the function  $g $ can be chosen as follows
$$
g(x)={\rm Ind}_{\Zp}(x)+ \frac{1}{ \delta_s^\phi } \sum_{n=1}^{\infty} \abs{s}^{\nu n}
\frac{1}{\chi({\rm S_n})}{\rm Ind}_{\rm S_n}(x),
$$
where $ {\rm Ind} _B $ denotes the indicator of the set $B $.
To conclude the proof, it remains to note that
$\ImU_{O}(\Lambda)$  is an open set, since the holomorphic
map $\Ims^{\phi}$ is open (proposition 2.1 part I \cite{RS} ).
\Endproof

\newpage
\begin{figure}

\centerline{}
\psfile{16}{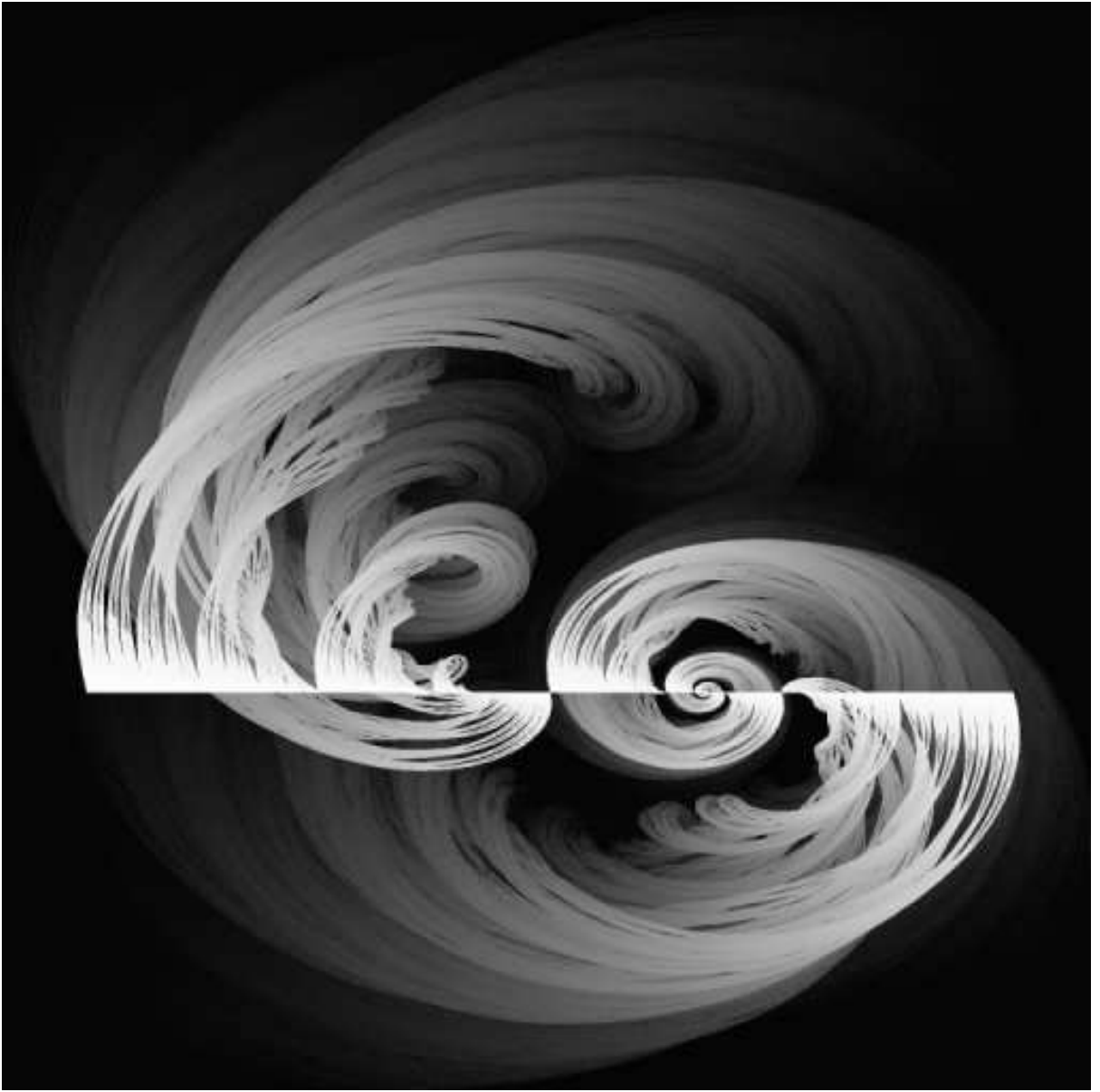}
\caption{
{ \footnotesize
The continuous transition between the 1- and 2-dimensional
spaces along the path drawn in figire 1
is shown.
}}\label{morph12}
\end{figure}
\begin{figure}
\psfile{16}{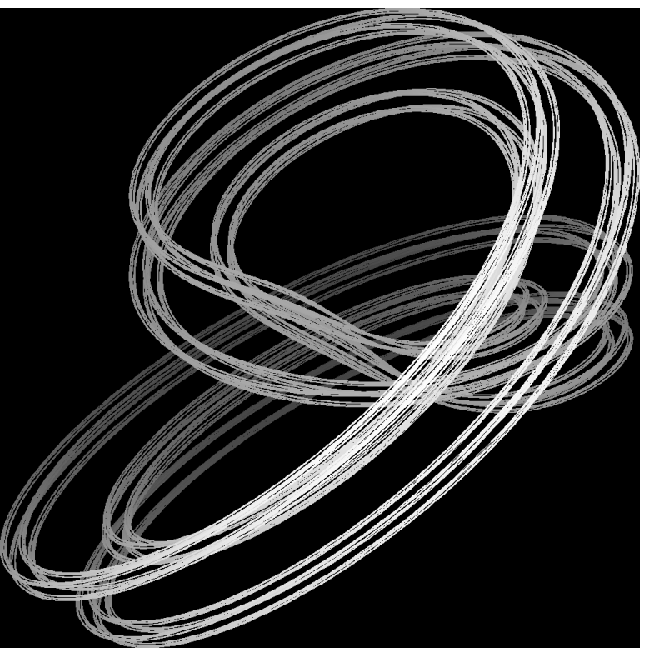}

\caption{
{\footnotesize
{  
The image of embedding  $ \Omega ^ {\phi} _ {\nu, \alpha}: \Sigma _ {\bm
{a}} \mapsto \Rs^3 $
at $a_k = (k+2)!$, $\nu=1.001$ ,$ \alpha=2$  and $\phi(t)=\exp(i2\pi t)$.
The Hausdorff dimension of this $\bm {a}$-adic solenoid image
$=1+\nu^{-1} \approx 1.999001$
}}\label{solenoid}
}

\end{figure}

\begin{thebibliography}{99}
\bibitem{Kob}
{\it N. Koblitz}\ p-Adic Numbers, p-Adic Analysis, and Zeta-Functions, Springer, New York Heidelberg Berlin (1977).
\bibitem{GG}
{\it I. M Gel'fand, M. I. Graev. and L. I. Pyatetskii-Shapiro }\
Representation Theory and Automorpttic Functions. Saunders. Philadelphia (1969).
\bibitem{VVZ}
{\it  V. S. Vladimirov, I. V. Volovich, and E. I. Zeienov}
\ p-Adic Analysis and Mathematical Physics, World Scientific, Singapore-New Jersey-London-Hong Kong (1994).
\bibitem{FED}
{\it H. Federer }\
Geometric Measure Theory, Springer, New York-Heidelberg-Berlin (1969).
\bibitem{HR}
{\it E. Hewitt and K. Ross}\
Abstract Harmonic Analysis, Vol. I, Springer, New York-Heidelberg Berlin (1963).
\bibitem{bill}
{\it P. Billingsley} \ Ergodic Theory and Information, Wiley, New York-London-Sidney (1965)
\bibitem{Feder}
{\it J. Feder} \ Fractals, Plenum, New York (1988).
\bibitem{Mam}
{\it D. Mumford.} \ Tata Lectures Notes on Theta Functions, Vols. I, II, Birkhauser, Boston-Basel-Stuttgart (1983. 1984).
\bibitem{Shabat}
{\it B. V. Shabat} \ Complex Analysis [in Russian], Vol. 1, Nauka, Moscow (1985).
\bibitem{Zelenov}
{\it Zelenov E.I. }//
J.Math.Phys. V32.147-152. 1991.
\bibitem{Pit}
{\it Pitkanen M.} // $p$-adic Physics. Department of Theoretical Physics,
University of Helsinki, SF-00170 Helsinki, Finland. 8. September 1994.
\bibitem{Speer}
{\it  Speer E. R. , Westwater M.J. } // Ann. Inst. Henri Poincare A14,1. 1971.
\bibitem{Hooft}
{\it  't Hooft G. , Veltman M. } // Nucl. Phys.  B44 189. 1972 .
\bibitem{Wilson}
{\it Wilson K. G. } // Phys. Rev. D7, 2911. 1973.
\bibitem{Mandel}
{\it Mandelbrot B.}\ General Property of Fractals.
Fractals in Physics, North Holland. Amsterdam Oxford-New York-Tokyo (1986).
\bibitem{MIS}
{\it Lerner E.U. , Missarov M. D.  } // Theoretical and Mathematical Physics.  v 101. N 2. 1994.
\bibitem{Chis}
{\it Chistyakov D.V. }// Theoretical and Mathematical Physics   Vol. 109. N 3. 1996
\bibitem{PBM}
{\it  Prudnikov A.P., Brichkov U.A., Marichev O.I. } \ Integrals and Powers  [in Russian], Nauka, Moscow ( 1981).
\bibitem{Bur1}
{\it Bourbaki  N.} \ Topologie Generale, Livre III .  Hermann.
\bibitem{Schwartz}
{\it Schwartz L.} \ Analyse Mathematique v.1.  Hermann (1967).
\bibitem{RS}
{\it Forster O.} \ Riemannsche Flachen. Springer-Verlag, Berlin Heidelberg NewYork (1977).
\end{thebibliography}
\end{document}